  \documentclass{article}
\usepackage{amssymb,amsmath,mathrsfs,amsthm,eepic,hyperref}
\usepackage{subfigure,graphicx,epsfig,color}

\newtheorem{theorem}{Theorem}[section]

\newtheorem{lemma}{Lemma}[section]
\newtheorem{rem}{Remark}[section]

\title{
         Homogenization Model for Aberrant Crypt Foci
\footnote{This work was partially supported by the project
PTDC/MATNAN/0593/2012
 and also 
CMUC -- UID/MAT/00324/2013, funded by the Portuguese Government through FCT/MCTES and co-funded by the European Regional Development Fund through the Partnership Agreement PT2020.
}}

 \author{ 
 {\small Isabel  N.~Figueiredo}\footnote{CMUC, Department of Mathematics, University
of Coimbra, Portugal (e-mails: isabelf@mat.uc.pt, carlosl@mat.uc.pt, roman@mat.uc.pt.).}
\and {\small Carlos Leal}\footnotemark[2]
\and {\small Giuseppe Romanazzi}\footnotemark[2]
\and {\small Bjorn Engquist}\footnote{Department of Mathematics and the Institute for Computational Engineering and Sciences, University of Texas at Austin, USA (email: engquist@ices.utexas.edu).} }
\date{\today}

\begin{document}

\maketitle{}

\noindent \textbf{Abstract--}{
Several explanations can be found in the literature about the origin of
colorectal cancer. There is however some agreement on the
fact that the carcinogenic process is a result of several genetic
mutations of normal cells. The colon epithelium is characterized by millions of invaginations,
very small cavities, called crypts,  where most of the cellular
activity occurs.
It is consensual in the medical community, that  a potential  first manifestation  of
the carcinogenic process, observed in conventional colonoscopy images, is the appearance
of Aberrant Crypt Foci (ACF). These are clusters of abnormal crypts,  morphologically characterized by  an atypical behavior of the cells that populate the crypts. In this work an homogenization model is proposed,  for  representing the cellular dynamics in the colon epithelium. The goal is  to simulate and predict, {\it in silico}, the spread and evolution of ACF, as it can be observed in colonoscopy images.
By assuming that the colon is an heterogeneous media,  exhibiting a periodic distribution of  crypts,  we start this work by describing a periodic model,  
that represents the ACF  cell-dynamics in a two-dimensional setting.
Then, homogenization techniques are applied  to this  periodic model,
  to find a simpler model,  whose solution symbolizes the averaged behavior of ACF at the tissue level.
  Some theoretical results concerning the existence of solution of the homogenized model are proven, applying a fixed point theorem. Numerical
results showing  the convergence of the periodic model to the
homogenized model are presented.
} \\

\noindent \textbf{Keywords - }{Homogenization,
convection-diffusion equations,  fixed-point theorem,  cell dynamics, colon.}

\smallskip


 \medskip



 \section{Introduction}\label{sec:intro}

 Colorectal cancer is one of the most common types of malignant
tumors in the Western World \cite{SW}. It is generally accepted that
its origin is associated with an accumulation of genetic mutations
at the cellular level that occur inside small cavities, called crypts, located
in the colon epithelium.

Despite the high rate of cancer mortality after detection, it is
possible to reduce the effects of this disease through an early
diagnosis.  This is due to the amount
of time,  20 to 40 years according to \cite{MILN}, that elapses between
the genetic mutations, that are at the origin of this process, and
the outbreak of the carcinoma. During this process some adenomas
(that are benign epithelial tumors)  develop  and if
detected and removed cancer can be avoided.
 
There are different ways to address the morphogenesis of the colorectal cancer.
 Most  of the approaches rely on experimental works (see {\it e.g.} \cite{Bakeretal})
that are then traduced in mathematical models.
The reader can refer to
 \cite{HJ,Kershaw,MILN,vLBJK,vLEIB,walter} 
for
a review of colorectal cancer models. It is believed that the
precursors of colorectal cancer are  Aberrant Crypt Foci (ACF)
\cite{birdI, birdII}. These latter are  clusters  of  crypts   in the colon epithelium, containing cells with a deviant behavior with respect to the normal ones.
There is no agreement regarding the morphogenesis of ACF.
 Two
 major  mechanisms have been proposed:
 top-down and  bottom-up  morphogenesis.
In the first, see \cite{lamlip,Sal}, the abnormal cells are, at the beginning, in a superficial portion of the mucosae and
after spread laterally and downwards inside the crypt. The bottom-up
morphogenesis \cite{lamlip,Pretal} relies on the property that the bottom of
the crypt is the most active region for cell division and genetic alterations and therefore assumes that  the appearance of   abnormal cells occurs in this region.
In \cite{Pretal} a combination of these two
mechanisms has been considered,  wherein an abnormal cell located in the crypt
base migrates to the crypt orifice, as in
the bottom-up morphogenesis,  and then it spreads
laterally filling the adjacent crypts in agreement to the top-down
morphogenesis. For a review of  ACF medical analysis and
colorectal cancer the reader can refer to
\cite{figu2010,greavesal,hurlstoneal,Pretal,RMB,Sal,tayloral}.

\begin{figure*}[t!]
\centering
   {\includegraphics[width=6.5cm,height=5.0cm]{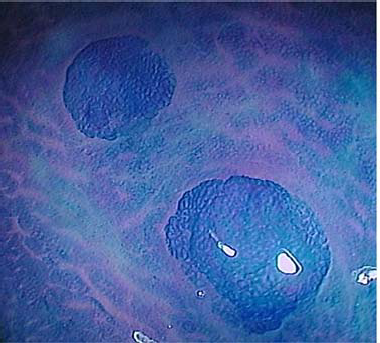}}
  \caption{Colonoscopy image (obtained with a conventional endoscope), showing two
 large  ACF in the human colon of a patient (Courtesy:  Faculty of Medicine, University of Coimbra, Portugal).}\label{fig:colon}
\end{figure*}

We  now briefly explain  our main motivation for doing this work and the main goal.
Conventional colonoscopy is a current medical technique used in gastroenterology to detect {\it in vivo}  ACF. The colonoscopy images give a top view of the colon wall, at the tissue level, and  from the endoscopic
point of view, ACF  stain darker  than normal crypts when some dyes (as for instance, methylene blue)  are instilled  in the colon (during the colonoscopy exam).
The medical Figure \ref{fig:colon}  shows a colonoscopy image exhibiting two ACF : the two
 large  dark round objects,
 each comprising  hundreds of aberrant crypts. In this figure we can also observe  the small circles spread all over the image that represent the top crypt orifices (directed to the lumen of the colon), and  we can clearly  perceive  a periodic distribution of crypts in the colon wall
(the dot like orifices that are densely distributed).
 
 The main motivation of this work is to model the evolution of ACF,  {\it i.e.}, the evolution of groups  of  several  aberrant  crypts, as observed in colonoscopy images. Our  goal is  to reproduce  {\it in silico} the dynamics of ACF at the macroscopic tissue level,  by using  appropriate mathematical models and tools.
To this end, we propose a cell dynamics model for describing  the evolution of abnormal colonic cells in a single crypt,
and subsequently, using the periodic structure of the colon, we resort to homogenization techniques for simulating
the evolution of ACF at the tissue level in a two-dimensional framework.
In particular, the cell dynamics model  used  is a PDE  (partial differential equation) model  that confines itself to two populations of colonic cells, normal and abnormal, and that is afterwards re-written only  in terms of the abnormal cell population.

 We give now a brief overview of the approach followed in this paper for modeling ACF dynamics. The colon is modeled as an heterogeneous two-dimensional
 medium  perforated by  crypts, that are periodically  distributed. The heterogeneous periodic model, used in this work, to represent the ACF
evolution in the colon is similar to that presented in  \cite{figu6}, 
 with the difference that in this work an hexagonal  region representing each crypt is used (as done in \cite{guebel2008computer}) instead of a square. The model is described 
 in detail in Sections  \ref{sec:geom}, \ref{sec:model} and \ref{sec:perio}. It consists of a system coupling a parabolic and an elliptic equation, whose unknowns are respectively, the density of abnormal cells and the pressure generated by   cell proliferation. We  assume that these abnormal cells have the property to proliferate not only inside the crypt as done by normal cells, but also outside. This property is believed to be responsible for abnormal cells to invade neighbor crypts and can induce crypt fission and the formation of an adenoma \cite{figu2009, Pretal, vLEIB}.
 In this work  crypt fission is not considered,  due to limitations of the chosen homogenization technique, a  two-scale asymptotic expansion method  (see Section \ref{sec:hom}).
 Then, by using an heuristic procedure, we derive the corresponding homogenized  model. It  is  again a system coupling a parabolic and an elliptic equation, whose unknown  is the  pair (abnormal cell density, pressure) corresponding to the  first terms  of the two-scale asymptotic expansions. Based on fixed-point type arguments,  the existence of solution of the homogenization model is proved. We remark that,  the complexity of the problem does not
 permit a proof of  the convergence of the periodic model to the homogenized one, though this is verified by the numerical tests we have performed.

After this introduction, the layout of the paper is as follows. In Section \ref{sec:geom} we define the three-dimensional (3D) crypt geometry as well as its reformulation
as a two-dimensional domain (2D), using a suitable bijection operator. In Section \ref{sec:model} the  crypt cell dynamics model   is described as well as its reformulation in a 2D framework. In Section \ref{sec:hom}  we extend this crypt cell dynamics model by periodicity  to obtain  an heterogeneous and periodic model for the colon. The homogenization technique, based on the heuristic  two-scale asymptotic expansion method ({\it ansatz}), is applied to this heterogeneous model to obtain the homogenization model.  In Section \ref{sec:solution},  the existence of solution to the homogenization model  is proved, as shown in Theorem \ref{tp}. The approximate solution to the homogenized model, that uses  finite elements for the space variable and finite differences for the time variable, is  described in Section \ref{sec:approx}.  The numerical simulations are presented in Section \ref{sec:simu}, and finally some conclusions and future work are outlined.

 \section{Crypt geometry (3D and 2D)}\label{sec:geom}

\begin{figure}[h!]
 \centering
  \includegraphics[scale = 0.47]{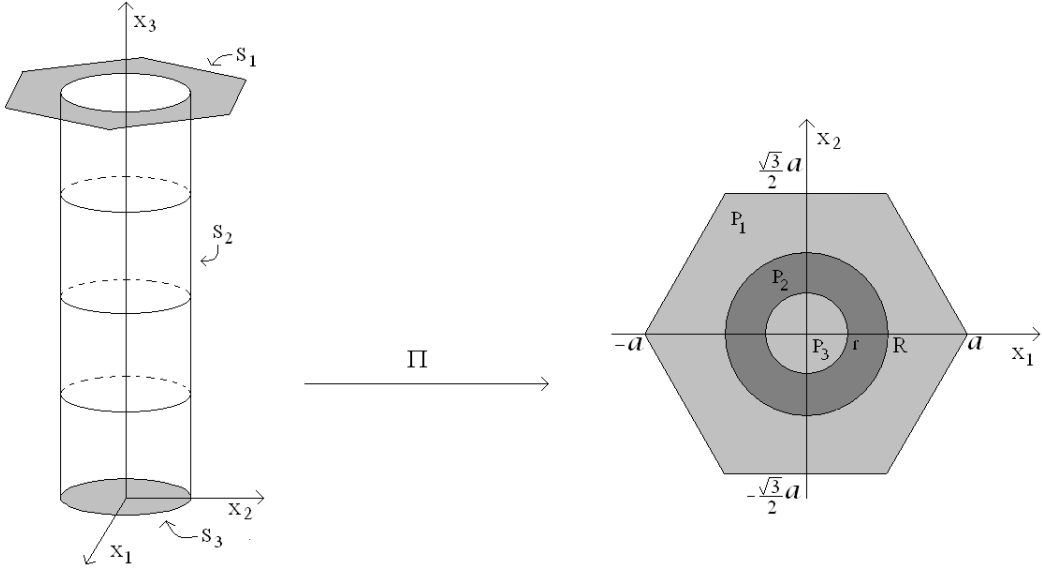}
 \caption{ Left: Schematic representation of a  region $S= S_1 \cup S_2 \cup S_3$ in the colon, containing  a 3D crypt and a small surrounding neighborhood (depicted by $S_1$). Right: The hexagonal region $P$ obtained by the projection of $S$ in a plane with the operator $\Pi$, where $P_i:= \Pi(S_i)$ for $i=1,2,3$.}\label{fig_geom}
 \end{figure}

Each crypt has a test-tube shaped structure, closed at the bottom and
 open  at the top.
According to  \cite{guebel2008computer, Halm},  the average dimensions for a human colonic crypt are 
 $433\, \mu m$ ($ 1 \mu m = 10^{-6}
  m$), from the bottom to the top, and $16\, \mu m$ for the diameter of the top orifice excluding the epithelium, with a cell  depth  of $15.1 \mu m$. In this work, we approximate the geometry  of a single crypt with a flat-bottomed cylinder $S\subset  \mathbb{R}^3$,  as  done in other works (see {\it e.g.} \cite{MFMB}). The cylinder $S$  is composed of three main parts $S:= S_1 \cup S_2 \cup S_3$ (represented in  Figure \ref{fig_geom}, left and defined in \eqref{1}). Note that we identify the crypt with the   region $S$,  where  $S_1$  (a regular hexagonal region of edge $a$) is the
inter-cryptal region around the crypt orifice, $S_2$ and $S_3$ represent,  respectively,  the lateral and the bottom region of the cylinder-like structure, with radius $R $
and height $L>0$ (for the numerical simulations, in Section \ref{sec:simu},  we assume a relation between the diameter $2R$ and the height $L$, in good agreement with the average  dimensions of a single crypt).
\begin{equation}\label{1}
\begin{array}{l}
S_1 := \left\{(x_1,x_2,x_3) \in \mathbb{R}^3 : (x_1,x_2) \in
{\rm regular\ hexagon\ of\ edge\ } a ,\ \sqrt{x_1^2+x_2^2}> R,\
\mbox{and}\ x_3 = L \right\}
\\
[.75em]
S_2 := \left\{(x_1,x_2,x_3) \in \mathbb{R}^3 : \sqrt{x_1^2+x_2^2} = R
\,\,\,\mbox{and} \,\,\,x_3 \in ]0,L]  \right\}
\\
[.75em]
S_3 := \left\{(x_1,x_2,x_3) \in \mathbb{R}^3 : \sqrt{x_1^2+x_2^2}
\leq R \,\,\,\mbox{and} \,\,\,x_3 =0  \right\}.
\end{array}
\end{equation}

For the subsequent work  herein (in order to define a 2D model for the ACF evolution)  it is necessary to reformulate  the  crypt geometry $S\subset  \mathbb{R}^3$ as a plane domain  
 $P$, that is defined as a regular hexagon, with area, $|P|$,  equal to 1, and with edge $a=\frac{\sqrt{2}}{3^{\frac{3}{4}}}$,
by using the  following bijective projection operator   $\Pi$ from
$S$ to $P $, defined by

\begin{equation}\label{2}
\Pi(x_1,x_2,x_3) :=
\begin{array}{lccl}
  \left\{
\begin{array}{lcl}
(x_1,x_2)  &\rm if & (x_1,x_2,x_3) \in S_1\\
[.75em]
\displaystyle\left(\frac{1}{R}\left[r+\frac{x_3}{L}(R-r)\right]x_1,\frac{1}{R}\left[r+\frac{x_3}{L}(R-r)\right]x_2\right)  & \rm if & (x_1,x_2,x_3) \in S_2\\
[.95em]
\displaystyle\frac{r}{R}(x_1,x_2)  &\rm if  & (x_1,x_2,x_3) \in S_3.
\end{array}
\right.
\end{array}
\end{equation}
The  projection $\Pi$  applies $S$ in  $P:= P_1 \cup P_2 \cup P_3$, where $P_i=\Pi(S_i)$.
We remark that $\Pi$ shrinks the bottom $S_3$ of the crypt into a circle with radius $r$ (verifying $0< r< R$), defined by $P_3$,
 while the lateral part $S_2$ of the crypt is projected in the region $P_2$, that lies  between the circles with radius $r$ and $R$.

 For any arbitrary function $g$ defined in $S\subset  \mathbb{R}^3$,   the corresponding function $g^*$ defined in \break  $\Pi(S)=P\subset  \mathbb{R}^2$ is given by
\begin{equation}\label{cv}
 g^*(X_1,X_2) = g(x_1,x_2,x_3) \quad {\rm with}  \quad (X_1,X_2)= {\Pi}(x_1,x_2,x_3).
 \end{equation}
Moreover the following relations, between the old and new variables, hold
\begin{equation}
\begin{array}{lcl}
\mbox{in} \quad  S_1 &  & (x_1,x_2)=  (X_1,X_2)  , \\
  \\
  \mbox{in} \quad  S_2 & &
\left \{\begin{array}{l}
  x_1 =   \displaystyle \frac{R \, X_1}{ \sqrt{X_1^2+ X_2^2}  } \\
\\
  x_2 =   \displaystyle \frac{R \, X_2}{ \sqrt{X_1^2+ X_2^2}  } \\
\\
 x_3 =  \displaystyle  \frac{L}{R-r} \big  ( \sqrt{X_1^2+ X_2^2} -r \big ), \\
\end{array}  \right.\\
\\
 \mbox{in} \quad  S_3 &  & (x_1,x_2)= \displaystyle \frac{R}{r} (X_1,X_2) .
\end{array}
\end{equation}
As a consequence,   the following relations between the partial derivatives of $g$ and $g^*$ also hold
\begin{equation}\label{der}
\begin{array}{lcl}
\mbox{in} \quad  S_1 & &
 \displaystyle  \big (\frac{\partial g}{\partial x_1},\frac{\partial g}{\partial x_2}, \frac{\partial g}{\partial x_3} \big  )=
  \displaystyle \big  (\frac{\partial g^*}{\partial X_1},\frac{\partial g^*}{\partial X_2}, 0 \big  ),   \\
  \\
  \mbox{in} \quad S_2 & &
\left \{\begin{array}{l}
  \displaystyle  \frac{\partial g}{\partial x_1} =   \displaystyle  \frac{1}{R}\left ( r+\frac{x_3}{L}(R-r)\right  )  \frac{\partial g^*}{\partial X_1} \\
\\
  \displaystyle  \frac{\partial g}{\partial x_2} =   \displaystyle   \frac{1}{R}\left ( r+\frac{x_3}{L}(R-r)\right  )  \frac{\partial g^*}{\partial X_2} \\
\\
  \displaystyle \frac{\partial g}{\partial x_3} =  \displaystyle  \frac{R-r}{RL} \big  ( x_1 \frac{\partial g^*}{\partial X_1}   +
   x_2 \frac{\partial g^*}{\partial X_2} \big ) ,\\
\end{array}\right. \\
\\
  \mbox{in} \quad S_3 & &
 \displaystyle(\frac{\partial g}{\partial x_1},\frac{\partial g}{\partial x_2}, \frac{\partial g}{\partial x_3})=      \frac{r}{R}  (\frac{\partial g^*}{\partial X_1},\frac{\partial g^*}{\partial X_2}, 0) .
\end{array}
\end{equation}



 \section{Colonic crypt cell model (3D and 2D)}\label{sec:model}

 Each crypt is a compartment containing different types of cells. These are aligned along the crypt wall: stem cells are believed to
reside in the bottom of the crypt, transit cells along the middle
part of the crypt axis and differentiated cells at the top of the
crypt.  In normal human colonic crypts,  the cells renew completely
each 5-6 days \cite{Ross}, through an harmonious and ordered procedure  which includes the
proliferation of cells, their migration along the crypt wall towards
the top and their apoptosis, as they reach the orifice of the crypt.

   It is generally believed that a sequence of genetical mutations and alterations in the  sub-cellular mechanisms
    can evolve later into colorectal cancer \cite{BC,MILN,vLBJK,walter}.
    In the literature there are many mathematical models proposed for
 reproducing the colonic cell dynamics.  We refer, for example, to the interesting and useful reviews   \cite{carulli2014unraveling,Kershaw,vLBJK}.  These models can be divided in two large groups: compartmental and spatial models.

 The compartmental models track the populations of different types (compartments) of colorectal cells (usually stem, semi-differentiated and differentiated  cells) by using systems of ordinary differential equations as in \cite{Boman2008,JEBMC}
or by stochastic models \cite{KomaWang,LBWP}.
These compartmental approaches determine, in each time, the number of cells of each population,  and allow cells to enter or leave the compartments without   specifying the relative spatial location of cells in the crypt. These models cannot describe then the cell migration and the spatial influence of proliferation or death rate parameters,  as opposed to spatial models that can incorporate these effects.
The group of spatial models can be  divided
 in two subgroups: cell-based models and continuum models.   In cell-based models, the cell dynamics is given by a finite number of mechanical forces and protein concentrations  among other quantities. Moreover, some probability of dividing or moving cells in certain directions can be used, as well as the  passage of  a mutation characteristic, therefore,  as in  \cite{FBC,Isele,MFMB,paulusal,Van2009}, stochastic models coupled with deterministic relations are  used.
Since these models treat cells as distinct entities they are able to incorporate sub-cellular features as cell-cell and cell-membrane adhesions, signaling pathways and protein level models \cite{Kershaw,walter}. The disadvantage of these approaches is that they are computationally expensive and need a lot of
experimental parameters that are difficult to
 collect.  Continuum models measure the cell populations by their densities and use   few parameters and computations to describe the continuum dynamics of cells.
Continuum models are usually defined by systems of partial differential equations to characterize the cell features and   are very popular
 for modeling avascular tumor growth \cite{APS,CLN,CLLW, RCM2007,SC, WK} or tissue growth \cite{KingFranks}. These models are used
 to simulate only the crypt cell dynamics \cite{figu4, MM, MWFETM,  Osborneetal,walter}
  or the crypt dynamics (budding and fission) \cite{DL2001,EC,figu5} based on a cell dynamics model.

 A continuum model was proposed in   \cite{figu6}  for describing the abnormal cell dynamics in the 3D  crypt domain $S\subset  \mathbb{R}^3$ and, through   its  reformulation in the 2D domain $P=\Pi(S) \subset  \mathbb{R}^2$ (these two domains were introduced in
 Section \ref{sec:geom}).
 In this section we define a similar model,  that will be used subsequently for deriving the homogenization
 model. 
 
  According to \cite{MILN} {\it ``Two scenarios of cellular dynamics within a colonic crypt are possible. In the first scenario, only colorectal stem cells are considered to be at risk of becoming cancer cells \cite{nowak2002role}. In the second scenario, all cells of the colonic crypt are assumed to be at risk of accumulating mutations that lead to cancer \cite{michor2004linear}. Mutations in specific genes might confer an increased probability to the cell to stick on top of the crypt instead of undergoing apoptosis''}. In this work, we assume this second scenario, since it   is  in good agreement with medical information obtained through the colonoscopy images and  correspondent biopsies, performed during the colonoscopy exam \cite{figu2009}.

 In healthy (normal) colonic crypts the cells are distributed in the following way: i) in the bottom there are stem cells with a high proliferative rate, ii) in the middle there are semi-differentiated cells, that also proliferate but with a lower rate than stem cells, iii) at the top of the crypt there are fully differentiated cells, that do not reproduce themselves.

 In our prototype model  we start by considering that there exist two main classes of colonic cells, normal  with density $N$
   and abnormal  with density $C$.
 A normal cell has a proliferative rate in good agreement with its location in the crypt (middle, bottom or top). An abnormal cell is a proliferative cell that  does not have a normal proliferative rate. For instance, if at the top of the crypt there are cells that reproduce themselves, then these are abnormal cells. 
 Furthermore we assume that the two types of cells verify the overall density hypothesis $N+C=1$ (and  $0\leq N, C\leq 1$).  This is equivalent to supposing  that no free-space exists and  that normal and
 abnormal cells 
 have the same volume.
   This is often used  in the context of living tissue growth \cite{KingFranks,Painter} and for modeling  tumor growth \cite{frie,RCS,SC,WK}.  This hypothesis is appropriate for modeling early stages of development of abnormal cells, as it is our intention in this work, but in advanced stages a  condition involving an increase of volume could also be envisaged.

   Let us denote  by $t$ the time variable belonging to the  interval   $[0,T]$, with $T>0$ fixed, and by   $N(x_1,x_2,x_3,t)$  and $C(x_1,x_2,x_3,t)$, respectively,   the normal and  abnormal cell densities,
 at  each point $ (x_1,x_2,x_3)$ of $S$ and at time $t$. Thus,  for the overall density hypothesis ,  $N(x_1,x_2,x_3,t)$  and $C(x_1,x_2,x_3,t)$ represent  the percentage of normal and abnormal cells at the spatial point $(x_1,x_2,x_3)$ and  time $t$.
Then, based  on models of tumor growth, described by systems of PDEs   and  relying on transport/diffusion/reaction models (see for instance
 \cite{KingFranks},  also \cite {frie,Osborneetal, RCM2007}  and \cite{BGM2008,CostaSole}),  we propose the following system of PDEs for representing the dynamics of these populations of colonic cells in  $ S\times(0,T) $
    \begin{equation}\label{eq:modelS2_S3}
  \left\{
  \begin{array}{lcl}
 \displaystyle\frac{\partial N}{\partial t} + \nabla  \cdot ( v_N \,
 N) & = & \nabla \cdot (D_N\nabla N)+     \gamma     N   ,\\
 [.75em]
 \displaystyle\frac{\partial C}{\partial t} + \nabla\cdot( v_C \, C)
 & = & \nabla \cdot (D_C\nabla C) +    \beta     C      , \\
  [.75em]
 N+C =1 .
 \end{array} \right.
 \end{equation}
\noindent Here $D_N,D_C$ are the diffusion coefficients of normal
and abnormal cells, respectively,  $ \gamma$  is the birth rate of
 normal cells and $ \beta $ the birth rate of abnormal cells,  that  depend on $x_3$, which is the crypt height ($\gamma =\gamma (x_3)$ and $\beta =\beta (x_3)$).
  The convective velocity of the normal and abnormal cells are denoted by $v_N$ and $v_C$, respectively. Finally $\nabla$ and $ \nabla  \cdot $ are   the gradient and divergence  operators, respectively.

  The  proliferative rate  $ \gamma$  of normal cells is known to decrease (see for instance \cite{DL2001}) with respect the height position in the crypt (along its vertical axis) and it is null in the upper third part of the crypt along its axis (see \eqref{values} for the definition of $\gamma$).
  The proliferation rate $\beta$ of abnormal cells is assumed to be higher than the proliferative rate  $ \gamma$ of normals cells  {($0\leq \gamma\leq \beta$); in addition it is assumed  that abnormal cells can also proliferate near the top of the crypt  (see  \eqref{values}).

We note that in \eqref{eq:modelS2_S3} when no abnormal cells are present in the crypt, that is  $C=0$ in $S$,  the system reduces to the first equation. This equation  represents a convection-diffusion-reaction model for normal cells that can proliferate up to the two thirds of the crypt height, as assumed in the literature.
If we consider that normal and abnormal cells are coexisting,
 the normal dynamics of cells is changed  to include  abnormal cells.

We   suppose also that  the two populations of cells have the same convective velocity $v_N = v_C = v$ (as has been done  for instance in  \cite{RCM2007}),
and  that the interior of the colonic crypt  is ``fluid-like" and the cell dynamics
obey  to a Darcy's law
\cite{Greenspan,Osborneetal, RCS,WK}. Thus, the common convective velocity is  defined by
$ v   = - \mu \nabla p$,
where $p$ is an internal pressure and $\mu$ is a positive constant describing the viscous-like properties  of the medium. In the sequel we consider $\mu=1$, for sake of simplicity.

 We observe that \eqref{eq:modelS2_S3} is  a simplified model of the real biological phenomenon,  that intends to model the  very early stages of development of abnormal colonic cells. It is a model for the proliferation and movement of two different cell populations in a crypt, normal and abnormal, such that the total cell density is constant and therefore there are no gaps.  We  remark that this model  \eqref{eq:modelS2_S3}  belongs to the class of models  presented  in \cite[system (1.1)]{KingFranks}. This system exhibits
the conservation laws for two species of cells (malignant in the tumour context and normal), satisfying the no void condition  and having the same convective velocity.
 Our model (6) is precisely the  model of \cite[system (1.1)]{KingFranks},
with the  proliferative rates $ \beta $ for the abnormal cell density $C$ and $  \gamma $ for  the normal cell density $N$.    And similarly to our proposed model,  \cite[system (1.1)]{KingFranks}  is complemented by a constitutive law   \cite[equation (3.1)]{KingFranks}   that is precisely a Darcy's law.

Then, by summing the two  first  equations in  (\ref{eq:modelS2_S3}) and using the overall density hypothesis $N+C=1$
we   obtain from (\ref{eq:modelS2_S3})
 the following elliptic-parabolic coupled model in $ S\times(0,T]$,  whose unknown is the pair $(C,p)$:
\begin{equation}\label{eq:modelS2_S3_2}
 \left\{
 \begin{array}{l}
\displaystyle\frac{\partial C}{\partial t}-\nabla\cdot( \nabla p \, C)  =  \nabla
\cdot (D_C\nabla C)+    \beta  \,   C   ,\\
[.75em]
- \Delta p  = \nabla \cdot \big ((D_C -D_N )
\nabla C \big ) +   (\beta - \gamma) \, C + \gamma  ,   
\end{array}
  \right.
\end{equation}
where $\Delta$ is the Laplace operator.
We denote by $p(x_1,x_2,x_3,t)$   the pressure  generated by
 cell proliferation,
 at each point $ (x_1,x_2,x_3)$ of $S$ and at time $t$.
After
introducing  the following new parameters $D$ and  $E$  defined by
\begin{equation}\label{D1D2}
\begin{array}{lcl}
D:= D_C,
  &\quad &
  E:=
 D_C-D_N  ,
\end{array}
\end{equation}
 the system  \eqref{eq:modelS2_S3_2} becomes
 \begin{equation}\label{eq:3dsystem}
\left\{
 \begin{array}{l}
\displaystyle\frac{\partial C}{\partial t}-\nabla \cdot( \nabla p \, C) = \nabla
\cdot (D\, \nabla C)+ \beta  \,   C  \\
[.5em]
 - \Delta p = \nabla\cdot (E\,  \nabla C) +   (\beta - \gamma) \, C + \gamma  . 
\end{array}
\right.
\end{equation}
The coefficients $\gamma $ and $\beta $ are assumed to be regular enough functions defined in the spatial domain  $S $. This regularity hypothesis is needed in Section \ref{sec:solution} for proving the existence of solution of the homogenization model.
The definitions and values for these parameters are given in Section \ref{sec:simu}.

Furthermore,  in Section \ref{sec:asymp}, for simplifying the computations in the derivation of the homogenization model, we assume  $D_N =D_C$ and consequently $E=0$ (see \eqref{constantDs}). Therefore,  we have that  at the top of the crypt, in $S_1$ where $\gamma=0$, the system \eqref{eq:3dsystem} becomes
 \begin{equation}\label{eq:3dsystemS1}
\left\{
 \begin{array}{l}
\displaystyle\frac{\partial C}{\partial t}-\nabla \cdot( \nabla p \,  C) = \nabla
\cdot (D \nabla C)+  \beta  \,   C   \\
[.5em] - \Delta p =  \beta  \,   C.
\end{array}
\right.
\end{equation}
meaning  that in $S_1$   the abnormal cells can proliferate (due to the presence of the  proliferation rate $\beta $ that is assumed to be  non-zero, see \eqref{values}) and that they are under the influence of a  convective term,  if    $  \nabla p  $  is non-zero.  

 It remains to reformulate  the model \eqref{eq:3dsystem}  in  the 2D domain $\Pi(S) =P \subset  \mathbb{R}^2$ introduced in Section \ref{sec:geom}. Using the projection $\Pi$ and
formulae \eqref{der}, relating the space derivatives $\frac{\partial g}{\partial
x_i}$ and $\frac{\partial g^*}{\partial X_j}$, the coupled parabolic-elliptic system   \eqref{eq:3dsystem} can be rewritten  in the 2D space-time  domain $P \times ]0,T]$ as follows:

\begin{equation}\label{eq:2dsystem}
\left \{\begin{array}{l} \displaystyle \frac{\partial C^*}{\partial
t}- \mathcal{A}^*_{ij}\frac{\partial}{\partial X_i}(C^*
\frac{\partial p^*}{\partial
X_j})=\mathcal{A}^*_{ij}\frac{\partial}{\partial X_i}(D^*
\frac{\partial C^*}{\partial X_j})+ 
\beta^* \,   C^*
\\  
 \\
\displaystyle -\mathcal{A}^*_{ij}\frac{\partial ^2 p^*}{\partial X_i
\partial X_j}=\mathcal{A}^*_{ij}\frac{\partial}{\partial X_i}(E^*
\frac{\partial C^*}{\partial X_j})  +
 (\beta^* - \gamma^*)  \,   C^*+
\gamma ^*
. \\  
\end{array}
\right.
\end{equation}
Here  $D^*$,  $E^*$, $\beta^*$, $\gamma^*$ and $D$,  $E$, $\beta$, $\gamma$  are related by   \eqref{cv} and
 the parameters $\mathcal{A}^*_{ij}$ for $i,j=1,2$   are defined in $P$ by
\begin{equation}\label{9}
\begin{array}{ll}
 \mathcal{A}^*_{ii}(X_1,X_2):=
\left \{
\begin{array}{lcl}
1 & \mbox{in} & P_1\\
[.5em]
g_i(X_1,X_2) & \mbox{in} & P_2\\
[.5em]
\left(\frac{r}{R}\right)^2 & \mbox{in} & P_3, \\
 \end{array}
 \right.
 &\quad
 \mathcal{A}^*_{12}(X_1,X_2) :=
\left \{
\begin{array}{lcl}
0 & \mbox{in} & P_1\\
[.5em]
g_3(X_1,X_2) & \mbox{in} & P_2 \\
[.5em]
0 & \mbox{in} & P_3,
\end{array}
\right.\\
\end{array}
\end{equation}
and $\mathcal{A}^*_{21}(X_1,X_2)= \mathcal{A}^*_{12}(X_1,X_2)$  with
\begin{equation}\label{10}
\begin{array}{l}
g_1(X_1,X_2) := \displaystyle \frac{X_1^2 + X_2^2}{R^2}+
\frac{(R-r)^2}{L^2(X_1^2
+ X_2^2)}X_1^2\\
[.95em]
g_2(X_1,X_2):= \displaystyle  \frac{X_1^2 + X_2^2}{R^2}+
\frac{(R-r)^2}{L^2(X_1^2
+ X_2^2)}X_2^2\\
[.95em]
g_3(X_1,X_2):= \displaystyle \frac{(R-r)^2}{L^2(X_1^2 + X_2^2)}X_1 X_2.
\end{array}
\end{equation}
We remark that the coefficients $A^*_{ij}$ satisfy the following uniform ellipticity condition
\begin{equation}\label{elli}
\mathcal{A}_{ij}^*(X_1,X_2) \xi_i \xi_j \geq \eta \xi_i
\xi_i,\,\,\,\,\,\, \forall(X_1,X_2) \in P \,\,\, \mbox{and}\,\,\,\,
(\xi_1,\xi_2) \in \mathbb{R}^2,
\end{equation}
with $\eta>0$ a constant.
In fact,  it is clear that  \eqref{elli} is verified in $P_1 \cup P_3$. In $P_2$,
 the inequality  $r^2 \leq X_1^2 + X_2^2$  yields
$$\
\begin{array}{lcl}
 \mathcal{A}_{ij}^*(X_1,X_2) \xi_i \xi_j &=& \displaystyle  \frac{X_1^2 + X_2^2}{R^2}(\xi_1^2+\xi_2^2)+
\left(\frac{(R-r)}{L \sqrt{X_1^2 + X_2^2}}X_1\xi_1 +\frac{(R-r)}{L
\sqrt{X_1^2 + X_2^2}}X_2\xi_2\right)^2 \\
&\geq &
 \displaystyle  \left(\frac{r}{R}\right)^2(\xi_1^2+\xi_2^2).
\end{array}$$
 This property \eqref{elli} guarantees the existence and the
 uniqueness  of the solution $p$ of the elliptic problem in the second equation of system \eqref{eq:2dsystem}, when considering its right hand side known.

 \section{Homogenization model}\label{sec:hom}

 Based on  strong biological and medical evidence, the morphology of the colon epithelium is characterized by millions of crypts (according to \cite{MILN} there are approximately $10^7$ crypts in the mammalian colon epithelium).  Therefore, it can be defined  mathematically   as an heterogeneous material with a periodic distribution of small heterogeneities: the crypts. The size of each crypt is measured in micro-meters,  and is very small in comparison with the size of the colon, that can be measured in centimeters.

The main goal of this work is to simulate the dynamics of abnormal cells in a region  of the colon epithelium, and somewhat describe  and predict the spread and evolution of ACF.

 Due to the periodic feature of the colon epithelium, mentioned above, we propose to extend by periodicity the abnormal cell dynamics  model \eqref{eq:2dsystem}  to find an heterogeneous model for the colon. Subsequently, and for analysing the macro-behavior  of the colon epithelium induced by its crypt micro-properties,
 we apply an homogenization technique.
  The homogenization method  used  is a   two-scale asymptotic expansion,   
suited for  models  posed in a periodic domain. Because of this  periodicity constraint we do not consider the possibility of crypt fission, that is advocated, for instance, in the ``bottom-up'' theory  for the morphogenesis of ACF (see \cite{greavesal,Pretal}). 
  With the adopted homogenization  technique
 we obtain a new model  (simpler and easier to handle),  called
 the homogenization model, that leads to an average (or equivalently,  homogenized) fictitious colon with the  corresponding homogenized abnormal cell behavior.  We refer for instance to \cite{bakhvalov1989homogenisation, bensoussanasymptotic,cioranescu1999introduction,Eng, sancheznon,tartar2009general} for a detailed presentation of the homogenization theory.

 We emphasize that with  the homogenization technique  we can perceive alterations at the surface of the colon epithelium, by   using the information of what happens inside the crypts at the cellular level. In particular,  in \cite[p. 202]{RCM2007},    it is commented on the relevance of using homogenization techniques to address problems related to the passage from discrete to continuous levels (in the present framework the discrete level corresponds to  the crypts and the continuous level to the tissue level, that  is the colon wall).

  \subsection{Heterogeneous periodic model}\label{sec:perio}

In this and in the following sections,  a region (piece) of the colon epithelium, hereafter denoted by
$\Omega \subset  \mathbb{R}^2$  
is represented by the heterogenous domain, obtained by the periodic distribution with basis $\varepsilon (\frac{3}{2} a, \frac{\sqrt{3}}{2}a)$ and $\varepsilon(0,\sqrt{3} a )$,
of a rescaled crypt 
 $\varepsilon P$.   Here $\varepsilon P$
 is a regular hexagon with edge of size $ \varepsilon a$
where $P$ stands for the  2D crypt  hexagon
 defined in Section \ref{sec:geom}. 
The parameter $ \varepsilon > 0 $ is very small
 when compared to the size of the domain $\Omega$. More precisely, we consider the surface of the colon epithelium as a 2D  periodic  structure,  since the domain $ \varepsilon P $,
  representing the {\it ``crypt  and the small region  surrounding its orifice"},  is replicated in all the domain $\Omega$, as depicted in Figure \ref{fig:crypts}, left.  We also emphasize that this  periodic  structure is observed in conventional colonoscopy images, that give a top view of the colon wall  at the tissue level, as exemplified in Figure \ref{fig:colon}, as well as in histological images    \cite{Araki1996, guebel2008computer}, were an hexagonal-like periodicity can be perceived, as displayed in Figure \ref{fig:crypts},  right.}

  \begin{figure*}[t!]
\centering
 {\includegraphics[width=6.0cm,height=4.0cm]{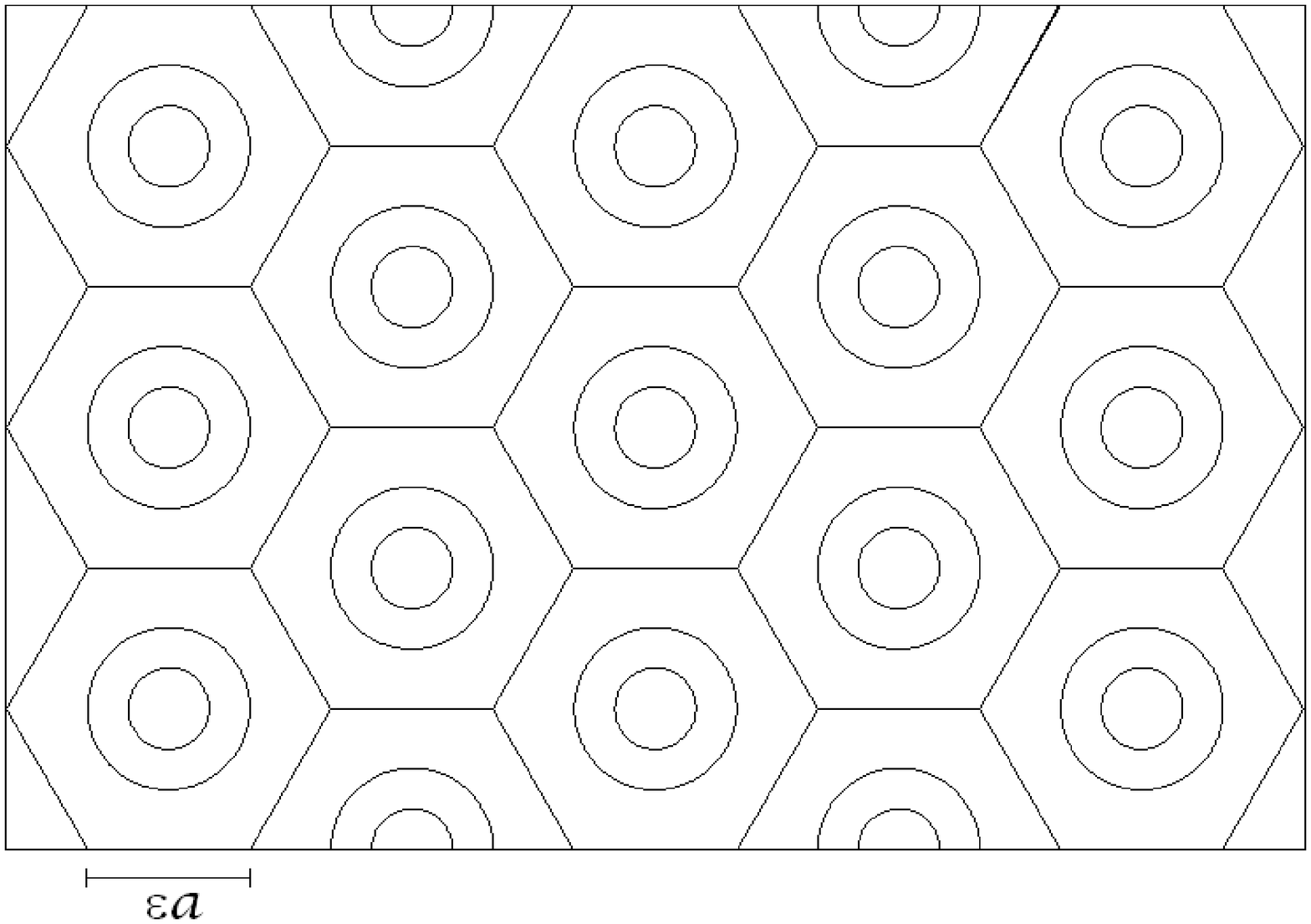}}\hspace{1cm}
{\includegraphics[width=6.0cm,height=4.0cm]{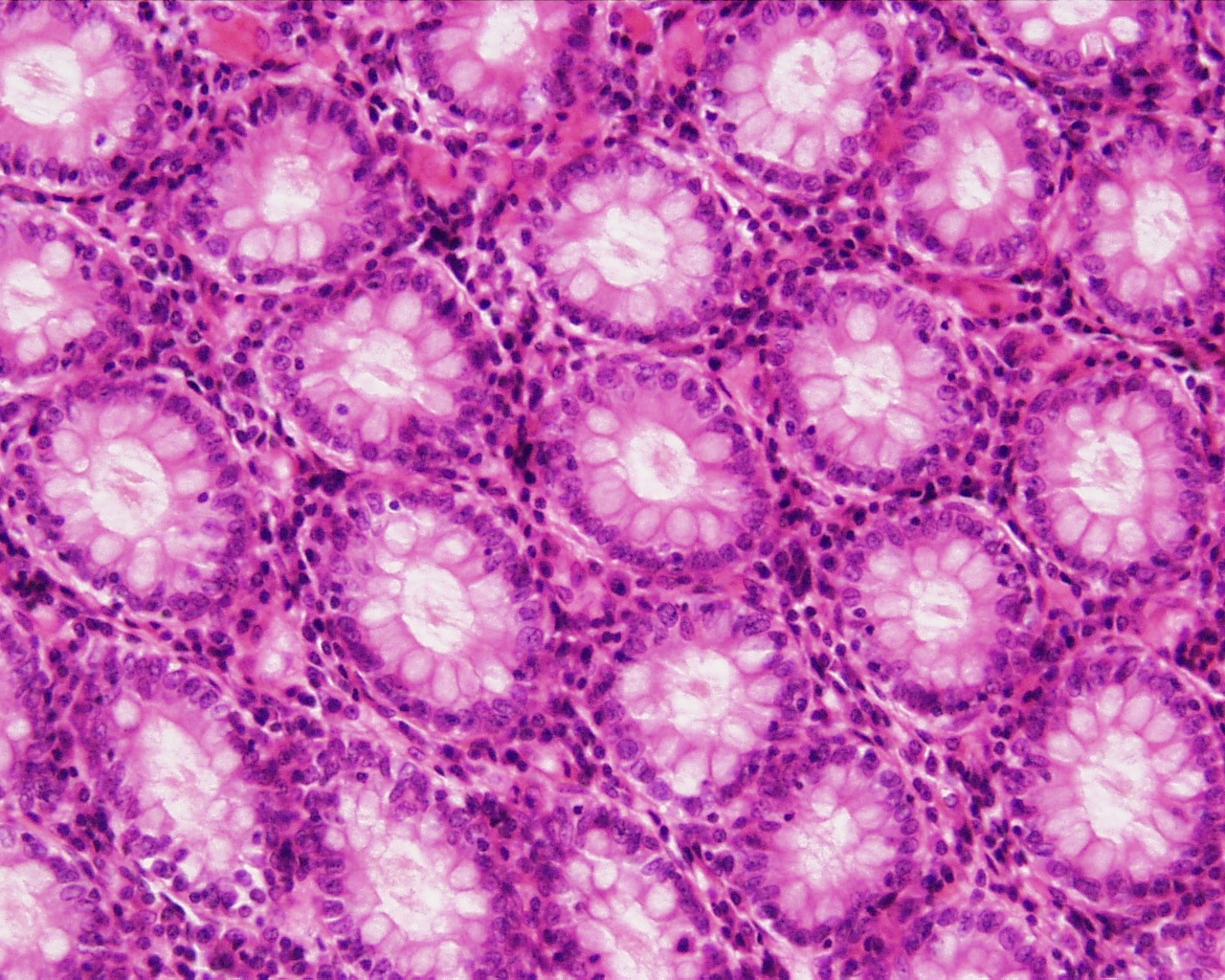}}
 \caption{ Left image: Two-dimensional and schematic  top view of part of the colon wall, 
 considered as an heterogeneous material with a periodic distribution of small heterogeneities,  $\varepsilon P$.
  Each  hexagon, with edge of size  $\varepsilon a$  
  represents the region  $ \varepsilon P $. Inside each hexagon,
  there is a crypt, containing the colonic cells, and that is represented by the two concentric circles;  the region inside  the hexagon  and outside the  large circle is the  region surrounding the crypt orifice
  Right image:  Histological image representing a transverse section of the human colon \cite{HillMA}.
 Compare $\varepsilon P$  with the single region $P$ shown in Figure \ref{fig_geom}, right.} \label{fig:crypts}
 \end{figure*}
Therefore we
define periodic coefficients  $A_{ij}^{\varepsilon},D^{\varepsilon},
E^{\varepsilon}, \gamma^{\varepsilon}, \beta^{\varepsilon}$ (based on the coefficients  introduced in \eqref{eq:2dsystem}) as follows:
 $$A_{ij}^{\varepsilon}(X)= A_{ij}\left(\frac{X}{\varepsilon}\right)\quad  {\rm and\ }\quad
A_{ij} = \left \{
\begin{array}{l}
\mathcal{A}^*_{ij}  \; \mbox{\ in\ } \; P,  \\
[.5em]
\mbox{extended by periodicity (with period $P$) }\\
\mbox{elsewhere in $\mathbb{R}^2$,}
\end{array}
\right.$$ where $X=(X_1,X_2) \in \Omega  \subset  \mathbb{R}^2$. Analogous
formulae are used for
defining the other coefficients $D^{\varepsilon}, E^{\varepsilon},
\gamma^{\varepsilon}, \beta^{\varepsilon}$.

Thus the   2D  heterogeneous  periodic model, representing the colon epithelium in a  2D  setting,  corresponds to the following
 system  of two partial differential equations,
of  elliptic and  parabolic type, defined in
$\Omega\times]0,T]$, where the unknown is the pair   cell density
$C^{\varepsilon}$ and pressure $p^{\varepsilon}$.
\begin{equation}\label{11}
\left \{
\begin{array}{l}
\displaystyle \frac{\partial C^{\varepsilon}}{\partial t}-
A_{ij}^{\varepsilon}\frac{\partial}{\partial X_i}(C^{\varepsilon}
\frac{\partial p^{\varepsilon}}{\partial
X_j})=A_{ij}^{\varepsilon}\frac{\partial}{\partial
X_i}(D^{\varepsilon} \frac{\partial C^{\varepsilon}}{\partial X_j})+
  \beta^{\varepsilon} \,
C^{\varepsilon}
 \\ 
 \\
\displaystyle -A_{ij}^{\varepsilon}\frac{\partial ^2 p^{\varepsilon}}{\partial X_i
\partial X_j}=A_{ij}^{\varepsilon}\frac{\partial}{\partial X_i}(E^{\varepsilon}
\frac{\partial C^{\varepsilon}}{\partial X_j})  +
  (\beta^{\varepsilon} -
\gamma^{\varepsilon}) \, C^{\varepsilon}+ \gamma ^{\varepsilon}
.  
\end{array}
\right.
\end{equation}
Homogeneous boundary conditions are assumed  for both
$C^{\varepsilon}$ and $p^{\varepsilon}$ in $\partial\Omega \times
]0,T]$, and the  initial condition, at time $t=0$,  for the  density of abnormal cells is represented by $C^{\varepsilon}(X,0) = C_0(X)$.
Moreover, compatibility conditions for $t=0$ are hypothesized, that is, $C_0(X)=0$
on $\partial \Omega $ and the pressure $p^{\varepsilon}(X,0) $  is the
solution to the elliptic equation in \eqref{11},   when
$C^{\varepsilon}$ is replaced by $C_0(X)$.

\subsection{Two-scale asymptotic expansions}\label{sec:asymp}

The homogenization technique applied to  \eqref{11}  considers a sequence of problems \eqref{11} indexed by $\varepsilon$,
and takes the limit of this sequence, when   $\varepsilon$ tends  to zero.  Equivalently, the objective is to find the limit  $(C^0,p^0)$ of the sequence of solution pairs $\{(C^{\varepsilon}, p^{\varepsilon})\}_{\varepsilon>0}$ of \eqref{11}, in an appropriate topological space. The homogenization model (or equivalently, the limit problem) is the model  which has solution  $(C^0,p^0)$.

Here, we use an heuristic procedure, that
consists of forming  a  two-scale asymptotic expansion in $\varepsilon$ for
$ C^{\varepsilon}$ and $ p^{\varepsilon} $  defined by
\begin{equation}\label{13}
C^{\varepsilon}(X,t)= \sum_{i=0}^{+\infty} \varepsilon^i C^i(X,Y,t) \quad
\mbox {and }\quad  p^{\varepsilon}(X,t)= \sum_{i=0}^{+\infty} \varepsilon^i p^i(X,Y,t) \qquad (\text {with} \; Y= X/ \varepsilon)
\end{equation}
in order to formally homogenize the system \eqref{11} (see \cite{bakhvalov1989homogenisation, bensoussanasymptotic,sancheznon}   for more details). 
In this section we identify the limit problem formally.

In \eqref{13}  $Y= X/ \varepsilon$   is  the microscopic variable, as opposed to $X$ that is the macroscopic  variable.
Each term $C^i(X,Y,t)$ and  $p^i(X,Y,t)$  are functions of the variables $X$,  $Y$, $t$,  and periodic in $Y$ with period $P$.

In what follows, and for the sake of  simplifying the computations we always assume that the
  diffusion coefficients are positive constants that have the same value for normal and abnormal cells.
   Since there is
very little  experimental quantitative information  about the value of these diffusion coefficients (see \cite{figueiredo2013physiologic} for a related work),
 it is common  in the literature to assume that these coefficients are constant and equal \cite{MILN, RCM2007}.
 So, from  \eqref{D1D2} we have
\begin{equation}\label{constantDs}
  D_N^{\varepsilon}(Y)=D_C^{\varepsilon}(Y)=D>0 \,\,\,\, \mbox{and}\,\,\,\, E^{\varepsilon}(Y)= 0,\ \forall Y\in  P ,
\end{equation}
and consequently the first term in the
second member of the second equation in system \eqref{11}
   disappears.

Then,  regarding $X$ and $Y$ as separate variables the differentiation operator $\frac{\partial}{\partial X}$  applied to a function depending on $(X,Y)$  becomes $\frac{\partial}{\partial X}+\frac{1}{\varepsilon}\frac{\partial}{\partial Y}$ and consequently
\begin{equation}\label{28_0}
\begin{array}{l}
 \displaystyle\frac{\partial p^{\varepsilon}}{\partial X_i} = \displaystyle
\frac{1}{\varepsilon}\frac{\partial p^{0}}{\partial Y_i}
+  \Big ( \frac{\partial p^{0}}{\partial X_i} + \frac{\partial p^1}{\partial Y_i} \Big )+
 \varepsilon   \Big ( \frac{\partial p^{1}}{\partial X_i} + \frac{\partial p^2}{\partial Y_i}   \Big ) +
\varepsilon^2  \Big ( \frac{\partial p^{2}}{\partial
X_i}  + \frac{\partial p^3}{\partial Y_i}    \Big ) + \mathcal O (\varepsilon^2) \\
\\
 \displaystyle\frac{\partial C^{\varepsilon}}{\partial X_i} = \displaystyle
\frac{1}{\varepsilon}\frac{\partial C^{0}}{\partial Y_i}
+\Big ( \frac{\partial C^{0}}{\partial X_i} + \frac{\partial C^1}{\partial Y_i} \Big )+
 \varepsilon   \Big ( \frac{\partial C^{1}}{\partial X_i} + \frac{\partial C^2}{\partial Y_i}   \Big ) +
\varepsilon^2  \Big ( \frac{\partial C^{2}}{\partial
X_i}  + \frac{\partial C^3}{\partial Y_i}    \Big ) + \mathcal O(\varepsilon^2)
\end{array}
\end{equation}
for $i=1,2$.  So, by introducing the two series   \eqref{13} into the equations   \eqref{11}  and identifying each coefficient
of   $\varepsilon^k$, for $k=-2,-1,0,1,\ldots$ we obtain  the
following  equations for the parabolic equation in  \eqref{11} :

\begin{itemize}

\item The $\varepsilon^{-2}$ equation

\begin{equation}\label{14}
- A_{ij}(Y)\left[\frac{\partial C^{0}}{\partial Y_i}\frac{\partial
p^{0}}{\partial Y_j}+C^0 \frac{\partial^2 p^{0}}{\partial Y_i
\partial Y_j}\right]=A_{ij}(Y)\left[ D \frac{\partial^2 C^{0}}{\partial Y_i
\partial Y_j}\right].
\end{equation}

\item The $\varepsilon^{-1}$ equation

\begin{equation}\label{15}
\begin{array}{rcl}
- \displaystyle A_{ij}(Y)\left[ \frac{\partial C^{0}}{\partial Y_i}
\left (\frac{\partial p^{0}}{\partial X_j}+\frac{\partial p^{1}}{\partial
Y_j} \right)+ \left  (\frac{\partial C^{0}}{\partial X_i}+\frac{\partial
C^{1}}{\partial Y_i} \right )\frac{\partial p^{0}}{\partial
Y_j}+C^1\frac{\partial^2 p^{0}}{\partial Y_i
\partial Y_j}\right. &  + & \\
& & \\
 \left.\displaystyle  C^0 \left( \frac{\partial^2 p^{0}}{\partial X_i
\partial Y_j}+\frac{\partial^2 p^{0}}{\partial Y_i
\partial X_j}+\frac{\partial^2 p^{1}}{\partial Y_i
\partial Y_j} \right)\right]& =& \\
& & \\
\displaystyle   D A_{ij}(Y)\left(\frac{\partial^2
C^{0}}{\partial X_i
\partial Y_j}+\frac{\partial^2 C^{0}}{\partial Y_i
\partial X_j}+\frac{\partial^2 C^{1}}{\partial Y_i
\partial Y_j}\right).  & &
\end{array}
\end{equation}

\item The  $\varepsilon^{0}$ equation

\begin{equation}\label{16}
\begin{array}{rll}
\displaystyle\frac{\partial C^0}{\partial t} -
A_{ij}(Y)\displaystyle
\Big [ \frac{\partial C^{0}}{\partial Y_i}
\left (\frac{\partial p^{1}}{\partial X_j}+ \frac{\partial p^{2}}{\partial
Y_j} \right ) + \left  (\frac{\partial C^{1}}{\partial X_i}+\frac{\partial
C^{2}}{\partial X_j} \right) \frac{\partial p^{0}}{\partial Y_j}
+ \left (\frac{\partial C^{0}}{\partial X_i}+\frac{\partial C^{1}}{\partial
Y_i} \right ) \left ( \frac{\partial p^{0}}{\partial X_j} +\frac{\partial
p^{1}}{\partial Y_j} \right )  &+& \\
& & \\
\displaystyle C^0\left(\frac{\partial^2 p^{0}}{\partial X_i
\partial X_j}+\frac{\partial^2 p^{1}}{\partial X_i
\partial Y_j}+\frac{\partial^2 p^{1}}{\partial Y_i
\partial X_j}+\frac{\partial^2 p^{2}}{\partial Y_i
\partial Y_j}\right)   & + &\\
& & \\
  \displaystyle  C^1 \left (\frac{\partial^2 p^{0}}{\partial X_i
\partial Y_j}+\frac{\partial^2 p^{0}}{\partial Y_i
\partial X_j}+\frac{\partial^2 p^{1}}{\partial Y_i
\partial Y_j} \right )+ C^2 \frac{\partial^2 p^{0}}{\partial Y_i
\partial Y_j}  \Big ]
& = & \\
& & \\
  D \displaystyle A_{ij}(Y)\left[ \frac{\partial^2 C^{0}}{\partial X_i
\partial X_j} +\frac{\partial^2 C^{1}}{\partial X_i
\partial Y_j}+\frac{\partial^2 C^{1}}{\partial Y_i
\partial X_j}
 + \frac{\partial^2 C^{2}}{\partial Y_i
\partial Y_j}
\right] +  \beta(Y) C^{0} . & &
\end{array}
\end{equation}

\end{itemize}

\noindent And  using the
same  procedure for the elliptic equation in  \eqref{11}  we obtain :

\begin{itemize}

\item The $\varepsilon^{-2}$ equation
\begin{equation}\label{17}
- A_{ij}(Y) \frac{\partial^2 p^{0}}{\partial Y_i \partial Y_j}=0.
\end{equation}

\item The $\varepsilon^{-1}$ equation
\begin{equation}\label{18}
- \displaystyle A_{ij}(Y)\left[ \frac{\partial^2 p^{0}}{\partial X_i
\partial Y_j}+\frac{\partial^2 p^{0}}{\partial Y_i
\partial X_j}+ \frac{\partial^2 p^{1}}{\partial Y_i
\partial Y_j}\right] =  0.
\end{equation}

\item The  $\varepsilon^{0}$ equation
\begin{equation}\label{19}
 -\displaystyle A_{ij}(Y)\left[\frac{\partial^2
p^{0}}{\partial X_i
\partial X_j}+\frac{\partial^2
p^{1}}{\partial X_i
\partial Y_j}+\frac{\partial^2
p^{1}}{\partial Y_i
\partial X_j}+\frac{\partial^2 p^{2}}{\partial Y_i
\partial Y_j}\right]   =
  (\beta(Y)-\gamma(Y)) C^{0}+ \gamma (Y).
\end{equation}
\end{itemize}

\noindent Since $p^{0}$ is assumed to be periodic in $Y$ and the coefficients  $A_{ij}^*$ verify the
uniform ellipticity condition \eqref{elli}, we have from \eqref{17} that $p^{0}(X,Y)=p^{0}(X)$.
So, from \eqref{18} we also obtain  $p^{1}(X,Y)=p^{1}(X)$. Equation
\eqref{14} gives $C^{0}(X,Y) = C^{0}(X)$ and consequently \eqref{15}
yields to  $C^{1}(X,Y) = C^{1}(X)$. Finally, these relations
permit us to rewrite  \eqref{19}   as the following  elliptic equation for $p^2$

\begin{equation}\label{20}
\displaystyle -A_{ij}(Y) \frac{\partial^2 p^{2}}{\partial Y_i
\partial Y_j}    =  A_{ij}(Y) \frac{\partial^2 p^{0}}{\partial X_i
\partial X_j}  +  (\beta(Y)-\gamma(Y)) C^{0}+ \gamma(Y)
.
\end{equation}

\noindent
Using the Fredholm alternative (see for instance
\cite{Eng}),  we can assert that   problem \eqref{20} has solution if,
for each $X\in\Omega$  
\begin{equation}\label{21}
\int_{P} [A_{ij}(Y) \frac{\partial^2 p^{0}}{\partial X_i
\partial X_j}(X)
+ (\beta(Y)-\gamma(Y)) C^{0}+ \gamma (Y)
]\, m\, dY =0
\end{equation}
where $m:\; P   \rightarrow \mathbb{R}$ is the unique solution of
\begin{equation}\label{22}
\left \{
\begin{array}{l}
-\displaystyle \frac{\partial^2 (A_{ij} m)}{\partial Y_i \partial Y_j}   =   0 ,
\\
 \\
\mbox{with} \; m \quad  P\mbox{-periodic},\;  m>0, \;  \displaystyle\int_P m \, dY = 1.
\end{array}
\right.
\end{equation}
Using the notation
   $\displaystyle\widetilde{f} = \frac{1}{|P|} \int_P f dY = \int_P f dY$ (because $|P|=1$), for
 the spatial average of a function $f$ in $P$, the
 condition \eqref{21} has the following form:

\begin{equation}\label{23}
 -\widetilde{A_{ij}m} \frac{\partial^2 p^{0}}{\partial X_i
\partial X_j}=   (\widetilde{\beta m} -\widetilde{\gamma  m} )C^{0}+ \widetilde{\gamma m}
   \quad  \mbox{in}\quad \Omega\times ]0,T].
\end{equation}

A  similar approach can be applied to equation \eqref{16}, as explained next.
Firstly, \eqref{19}  can be used to rewrite \eqref{16} as an elliptic equation for $C^2$ in the following way
 \begin{equation}\label{eq:C2elliptic}
\begin{array}{ll}
  \displaystyle D A_{ij}(Y) \frac{\partial^2 C^{2}}{\partial Y_i \partial Y_j} = &\\
  [.5em]
  \displaystyle \frac{\partial C^0}{\partial t} - A_{ij}(Y)
\left( \frac{\partial C^{0}}{\partial X_i} \frac{\partial
p^{0}}{\partial X_j} + D\frac{\partial^2 C^{0}}{\partial X_i\partial X_j} \right) - (\beta(Y)-\gamma(Y))C^0(1-C^0).
\end{array}
\end{equation}
Since $D$ is a positive constant and  using again the Fredholm alternative,  we have that  \eqref{eq:C2elliptic}
has a solution if,  for all $X\in\Omega$ the equation
\begin{equation}\label{eq:intpC0}
\int_P \left[\frac{\partial C^0}{\partial t} - A_{ij}(Y)
\left( \frac{\partial C^{0}}{\partial X_i} \frac{\partial
p^{0}}{\partial X_j} + D\frac{\partial^2 C^{0}}{\partial X_i\partial X_j} \right) - (\beta(Y)-\gamma(Y))C^0(1-C^0) \right]
mdY=0
\end{equation}
is verified, where $m$ is again the solution of \eqref{22}. Since  $C^0=C^0(X)$ and $p^0=p^0(X)$,  therefore \eqref{eq:intpC0} can be rewritten as follows
\begin{equation}
\frac{\partial C^0}{\partial t} - \widetilde{A_{ij}m}\frac{\partial C^{0}}{\partial X_i} \frac{\partial
p^{0}}{\partial X_j} =  D\widetilde{A_{ij}m}\frac{\partial^2 C^{0}}{\partial X_i\partial X_j} +
 (\widetilde{\beta m}-\widetilde{\gamma m}) C^0(1-C^0).
\end{equation}

In conclusion, the (formal) homogenization model   defined in $\Omega \times ]0, T]$, and resulting from the heterogeneous periodic model  \eqref{11},  is the following system
\begin{equation}\label{eq:homog}
\left \{
\begin{array}{l}
  \displaystyle \frac{\partial C^0}{\partial t} - \widetilde{A_{ij}m}\frac{\partial C^{0}}{\partial X_i} \frac{\partial
p^{0}}{\partial X_j} =  D\widetilde{A_{ij}m}\frac{\partial^2 C^{0}}{\partial X_i\partial X_j} +
 (\widetilde{\beta m}-\widetilde{\gamma m}) C^0(1-C^0)
\\
\\
 -  \displaystyle\widetilde{A_{ij}m} \frac{\partial^2 p^{0}}{\partial X_i
\partial X_j}=
  (\widetilde{\beta m} -\widetilde{\gamma  m} )C^{0}+
\widetilde{\gamma m} 
\end{array}
\right.
\end{equation}
whose unknown is the pair $(C^0,p^0)$. This is a
macroscopic model that represents the evolution
 of ACF at the surface of the colon, by using the  information of the  cell dynamics in the crypts. This information is transmitted to the homogenization model through the physiologic parameters  that  in \eqref{eq:homog}  are averaged (they are represented by the same letters as before with a tilde over them).

\section{Existence of solution to the homogenization model}\label{sec:solution}

\subsection{Some Banach spaces and corresponding norms}\label{bn}

We denote by $X_i$, with $i=1,2$,  the components of a generic point $X=(X_1,X_2) \in \Omega$
and by $ D_{X_{i}}$ the partial derivative with respect to $X_i$. For  $\alpha\in ]0,1[$  and $T>0$ fixed, we consider the Banach
spaces $ C^{\alpha}(\overline{\Omega})$ and   $ C^{\alpha}(\overline{\Omega}_T)$, where $\overline{\Omega}_T =
 \overline{\Omega} \times [0,T]$, endowed with the norms (see for instance \cite{Pao})
$$\| u\|_{C^{\alpha}(\overline{\Omega})} =
\sup_{X \in \overline{\Omega} }|u(X)| + \sup_{X,Y \in \overline{\Omega} }\frac{|u(X)-u(Y)|}{|X-Y|^{\alpha}}$$
and
$$\|u\|_{C^{\alpha}( \overline{\Omega}_T)} =
\sup_{(X,t) \in \overline{\Omega}_T}|u(X,t)| +
\sup_{(X,t),(\xi,\tau)\in \overline{\Omega}_T}
\displaystyle\frac{|u(X,t)-u(\xi,\tau)|}{(|t-\tau|+|X-\xi|^2)^{\frac{\alpha}{2}}}.$$
We denote by $H^{1+\alpha, \alpha}(\overline{\Omega}_T) $ the space
$$ H^{1+\alpha, \alpha}(\overline{\Omega}_T)= \left\{
u:\overline{\Omega}_T \longrightarrow \mathbb{R}\;  \Big | \quad u, \; D_{X_{i}}u \in C^{\alpha}(\overline{\Omega}_T), \;
  \| u \|_{H^{1+\alpha,
\alpha}(\overline{\Omega}_T)}< \infty   \right\}$$ equipped with the
norm
\begin{equation}\label{26}
\begin{array}{lcl}
\| u \|_{H^{1+\alpha, \alpha}(\overline{\Omega}_T)} & = &
\displaystyle \sup_{(X,t) \in \overline{\Omega}_T} |u(X,t)| \; +
\displaystyle \sup_{(X,t), (X,t') \in \overline{\Omega}_T}
\frac{|u(X,t)-u(X,t')|}{|t-t'|^{\frac{\alpha}{2}}}\\
 & & \\
 & + &\displaystyle \sum_{i=1,2}\left(
\sup_{ (X,t) \in \overline{\Omega}_T} |D_{X_{i}}u(X,t)|
+\sup_{(X,t),(Y,t) \in
\overline{\Omega}_T}\frac{|D_{X_{i}}u(X,t)-D_{X_{i}}u(Y,t)|}{|X-Y|^{\alpha}}  \right.\\
 & & \\
& + &\left.\displaystyle \sup_{(X,t),(X,t') \in
\overline{\Omega}_T}\frac{|D_{X_{i}}u(X,t)-D_{X_{i}}u(X,t')|}{|t-t'|^{\frac{\alpha}{2}}}
\right).
\end{array}
\end{equation}
We remark that $H^{1+\alpha, \alpha}(\overline{\Omega}_T)$ is a
Banach space (see for example \cite[p. 795]{fri1},  where
it is proven that $H^{1+\alpha, \alpha}(\overline{\Omega}_T)$ is a Banach space, when endowed with  a norm  equivalent to that defined in  \eqref{26}).

Moreover, if $V$ is  a Banach space,  we also consider the set
$$\mathcal{C}^{\alpha}([0,T];V) = \left\{u:[0,T]\longrightarrow V \;  \left| \quad  \frac{\| f(t+h)-f(t)\|_V}{|t-s|^{\alpha}}<
M,\quad \forall  t,t+h \in [0,T] \right.\right\}.$$
where $M$ is a constant. This space  is also a Banach space when equipped with the norm
\begin{equation}\label{27}
  \|u \|_{\mathcal{C}^{\alpha}([0,T];V)}= \sup_{t \in [0,T]}\|u(t)\|_{V} +
 \sup_{s,t \in [0,T]}\frac{\| u(t)-u(s)\|_{V}}{|t-s|^{\alpha}}.
 \end{equation}
In addition we also consider the following two spaces (see for example \cite{evans1998partial}):

\smallskip

\noindent  - the space $C^{1, \alpha}(\overline{\Omega})$ consisting of those functions $u: \Omega \rightarrow  \mathbb{R}$ that are one-time continuously differentiable and whose  first-partial derivatives are H\"older continuous with exponent $\alpha$.

\smallskip

\noindent  - the  Sobolev space $W^{2,q}(\Omega)$, for $q\geq 1$, consisting of all locally summable functions $u: \Omega \rightarrow  \mathbb{R}$,
whose weak partial derivatives, up to the order   $2$, exist and belong to $L^q(\Omega)$.

\subsection{Partial existence results}

The two following lemmas assure the existence and uniqueness of
solutions of the parabolic and elliptic  equations \eqref{eq:homog},
separately. These results  will be used in the next
section for proving the existence of solution of the homogenized system
\eqref{eq:homog}.

In what follows we will assume that the domain $\Omega$  is
sufficiently regular to be able to use the known results on
{\it a priori}
estimates for the solution of elliptic and semilinear parabolic
problems.

\begin{lemma}\label{lem-par}
Let  $W \in \left(C^{\alpha}(\overline{\Omega}_T)\right)^2$
 and consider the  semilinear parabolic problem

\begin{equation}\label{28}
\left \{
\begin{array}{lcl}
\displaystyle\frac{\partial C}{\partial t}  - LC=f(X,t,C), &
\mbox{in} & \Omega \times ]0,T] \\
& & \\
 C= 0,& \mbox{on} & \partial \Omega \times ]0,T]\\
& & \\
C(X,0) = C_0(X)& \mbox{in} & \Omega \times \{0\}
\end{array}
\right.
\end{equation}
where $L :=  \widetilde{A_{ij}m}  W \frac{\partial
}{\partial X_i}+ D\widetilde{A_{ij}m}\frac{\partial^2 }{\partial
X_i
\partial X_j}$ , $f(X,t,C) :=  (\widetilde{ \beta m}-\widetilde{\gamma m})C(1-C)$ and  the initial condition
$C_0 \in  C^{\alpha}(\overline{\Omega})$ satisfies
$C_0(X)_{|_{\partial\Omega}}= 0$ with $0 \leq C_0(X) \leq 1$. Then,
\begin{enumerate}
\item Problem \eqref{28} has  an unique solution $C$,  such that:
$$ 0\leq C(X,t) \leq 1 ,\quad  \forall (X,t) \in \overline{\Omega}_T.$$
\item
If the initial condition $C_0$ is regular enough and verifies
$C= 0$ on  $\partial \Omega$,   then $C \in H^{1+\alpha,\alpha}(\overline{\Omega}_T)$ and if
$(\widetilde{ \beta m}-\widetilde{\gamma m}) $ is small enough, there exists a positive constant $M$ such that
$$  \| C \|_{ H^{1+\alpha,\alpha}(\overline{\Omega}_T)} \leq M .$$
\end{enumerate}
\end{lemma}

\noindent \textbf{Proof:}  The first statement    is a direct
consequence of \cite[Theorem 5.1, p. 66]{Pao}.
Concerning the second statement, we  use   \cite[ Theorem 8, p. 204]{fri2}.  The hypotheses of this lemma for the domain,   initial
condition and  coefficients are the same  as those used in  \cite[ Theorem 8, p. 204]{fri2}, but for this latter theorem  there
is the following extra requirement   
for $f(.,.)$ in 
\eqref{28} that should be fulfilled :
\begin{equation}\label{fried_conition}
2 K | f(X,t,C )| \leq M  \quad  \mbox{in } \quad \Omega_T
\end{equation}
(where $|.|$ represents the modulus) for all $C \in H^{1+\alpha,\alpha}(\overline{\Omega}_T)$,
with $\|C\|_{  H^{1+\alpha,\alpha}(\overline{\Omega}_T)} <  M $,
where $K$ is a constant  depending on $\alpha$, $\Omega_T$ and on
the  coefficients of the convective term of the parabolic operator. Due  the particular form of the function $f$, in our case,  we
can conclude the condition \eqref{fried_conition} is satisfied if
$(\widetilde{ \beta m}-\widetilde{\gamma m}) $ is sufficiently
small. $\Box$

 \medskip

\begin{lemma}\label{lem-elli}
For a given  $ C\in C^{\alpha}(\overline{\Omega}_T)$ and for each
$t \in [0,T]$ the elliptic problem
\begin{equation}\label{29}
\left \{
\begin{array}{lcl}
\displaystyle -(\widetilde{A_{ij}m}) \frac{\partial^2 p}{\partial X_i
\partial X_j}=
 (\widetilde{\beta m} -\widetilde{\gamma  m} ) C +
\widetilde{\gamma m}  ,  & \mbox{in}
& \Omega \times [0,T]\\
& & \\
 p =
0,& \mbox{on} & \partial \Omega \times [0,T]
\end{array}
\right.
\end{equation}
has an unique solution $p$ in
$\mathcal{C}^{\frac{\alpha}{2}}([0,T];W^{2,q}(\Omega))$, for $q>2$.
Furthermore
\begin{equation}\label{30}
\|p\|_{\mathcal{C}^{\frac{\alpha}{2}}([0,T];W^{2,q}(\Omega))} \leq M
\end{equation}
where $M$ is a constant that depends on $\| C
\|_{C^{\alpha}(\overline{\Omega}_T)}$.
\end{lemma}

\noindent \textbf{Proof:} Using  regularity results and {\it a priori}
estimates for second order elliptic equations (see for instance
\cite[p. 95]{Pao}) and the regularity of $C$ we  conclude that, for
each fixed $t \in [0,T]$,  the solution $p$ of \eqref{29} is in
$W^{2,q}(\Omega)$ with $q>2$ and
\begin{equation}\label{elliptic-estimation}
\|p\|_{W^{2,q}(\Omega)} \leq M \| 
(\widetilde{\beta m} -\widetilde{\gamma  m} )C + \widetilde{\gamma m}
 \|_{L^p(\Omega)}
\end{equation}
where $M$ is a constant independent of $t$. Thus,  using  the regularity of $C$, we have
\begin{equation}\label{31}
\|p\|_{W^{2,q}(\Omega)} \leq M
\end{equation}
 where $M>0$ is another constant independent of
$t$ and depending on $C$, which can be small   if $ \|
\widetilde{\gamma m} \|_{L^p(\Omega)}$ and $\|
 \widetilde{\beta m} -\widetilde{\gamma  m}   \|_{L^p(\Omega)}$  are sufficiently small.

Considering now  \eqref{29}  for two different times
$t=t_1$ and $t=t_2$, subtracting the two corresponding equations and  dividing  by
$|t_1-t_2|^{\alpha}$,  we obtain
$$ -(\widetilde{A_{ij}m}) \frac{\partial^2 }{\partial X_i
\partial X_j}\left( \frac{p(t_1)-p(t_2)}{|t_1-t_2|^{\alpha}}\right)=
  (\widetilde{\beta m} -\widetilde{\gamma  m} )  
\left(\frac{C(X,t_1)-C(X,t_2)}{|t_1-t_2|^{\alpha}}\right),
\,\,\,\,\,\mbox{in} \,\,\,\, \Omega \times [0,T].$$ Then,   using
the same
 argument  as the beginning  of this proof, we have
\begin{equation}\label{32}
\left\| \frac{p(t_1)-p(t_2)}{|t_1-t_2|^{\alpha}}\right\|_{W^{2,q}(\Omega)} \leq
M
\end{equation}

\noindent where $M$ is another constant, depending on $\gamma$, but
independent of $t$. From \eqref{31}-\eqref{32} we have \break $p \in
\mathcal{C}^{\frac{\alpha}{2}}([0,T];W^{2,p}(\Omega))$ and also

\begin{equation}\label{33}
||p||_{\mathcal{C}^{\frac{\alpha}{2}}([0,T];W^{2,q}(\Omega))} \leq
M.  \qquad \Box
\end{equation}

\begin{rem}\label{rem2}
We observe that from \eqref{elliptic-estimation}-\eqref{33} we can conclude that, for
each $C \in {C^{\alpha}(\overline{\Omega}_T)}$, such that $ \| C
\|_{C^{\alpha}(\overline{\Omega}_T)} \leq K $, with $K$ constant, then the norm $\| p \|_{\mathcal{C} ^ {\frac{\alpha}{2}} ([0, T], W^{2,p}
(\Omega))} $ can be arbitrarily small  by taking $\widetilde{\gamma m} $ and   $\widetilde{\beta m} -\widetilde{\gamma  m}$ 
sufficiently small.
\end{rem}

\subsection{Main existence result}
Using  a fixed-point type argument we can prove the
 existence of a solution  to the homogenization model \eqref{eq:homog}.

\begin{theorem}\label{tp}
Let $\Omega$ be an open and bounded domain in $\mathbb{R}^2$ with
boundary $\partial \Omega $ regular enough,  $T$ a fixed positive
real number and   $(\widetilde{\gamma m}) $, $(\widetilde{\beta
m}-\widetilde{\gamma m})$ sufficiently small. Then,  the
homogenization model \eqref{eq:homog} has a solution $(C^0,p^0)$
such that $p^0 \in
C^{\frac{\alpha}{2}}([0,T];C^{1,\alpha}(\overline{\Omega}))$ and $
C^0 \in C^{\alpha}(\overline{\Omega}_T)$.

\end{theorem}

\noindent
\textbf{Proof}: The proof  is based on Schauder's fixed
point theorem (see for instance \cite{evans1998partial}). So the goal is to define an appropriate continuous operator $\mathcal {F}$    in a
compact and convex  subset $K$  of a Banach space, such that $\mathcal {F}: K\rightarrow K$, which implies that $\mathcal {F}$ has a fixed point.

Thus we first consider an operator $S$ defined in the domain
\begin{equation}\label{34}
B_M(0) = \{ v \in C^{\alpha}(\overline{\Omega}_T) : v(X,0)= C_0(X)
\,\,\, \mbox{and}\,\,\, || v ||_{C^{\alpha}(\overline{\Omega}_T)} \leq
M \}
\end{equation}
with $M$ a fixed constant and $C_0$ the initial condition in
\eqref{28}. We prove in the following that $\mathcal {F}: B_M(0)\rightarrow B_M(0)$, defined as $\mathcal {F}= Q \circ G \circ D \circ I
\circ E$,
verifies Schauder's fixed point theorem.

\begin{itemize}
\item   $E$ is the continuous operator defined by
\begin{equation}\label{35}
\begin{array}{lccl}
E: & B_M(0)
& \longrightarrow &\mathcal{C}^{\frac{\alpha}{2}}([0,T];W^{2,q}(\Omega)) \\
 & C & \longmapsto & p
\end{array}
\end{equation}
where   $q>2$. The image $p=E(C)$ is   the  solution of the elliptic problem \eqref{29}, given by Lemma \ref{lem-elli}.

\item   $I$ is an inclusion, that   identifies $p=E(C)
\in C^{\frac{\alpha}{2}}([0,T];W^{2,q}(\Omega))$ with $p \in
C^{\frac{\alpha}{2}}([0,T];C^{1, \alpha}(\overline{\Omega}))$. The
Sobolev embedding theorem  (see  for instance \cite[p. 270]{evans1998partial})
guarantees the continuity  of $I$.

\item $D$ is an operator that applies each  $ p\in \mathcal{C}^{\frac{\alpha}{2}}([0,T];C^{1,\alpha}(\overline{\Omega}))$
in $\hat{p} \in H^{1+\alpha,
\alpha}(\overline{\Omega}_T) $ defined by
\begin{equation}\label{36}
\begin{array}{lccl}
\hat{p}: & \overline{\Omega}_T
& \longrightarrow & \mathbb{R} \\
 & (X,t) & \longmapsto & \hat{p}(X,t) = p(t)(X).
\end{array}
\end{equation}
We can easily conclude, using the norms \eqref{26} and \eqref{27}
that there exist a constant $M>0$ such that
\begin{equation}\label{37}
\displaystyle ||\hat{p}||_{H^{1+\alpha,
\alpha}(\overline{\Omega}_T)}\leq M
||p||_{\mathcal{C}^{\frac{\alpha}{2}}([0,T];C^{1,\alpha}(\overline{\Omega}))}.
\end{equation}

\item $G$ is the gradient operator defined by
\begin{equation}
\begin{array}{lccl}
G : & H^{1+\alpha,
\alpha}(\overline{\Omega}_T)
& \longrightarrow & \left(C^{\alpha}(\overline{\Omega}_T)\right)^2  \\
[.5em]
 & \hat{p} & \longmapsto & \nabla\hat{p}= \displaystyle \left(\frac{\partial\hat{p}}{\partial X_1},\frac{\partial\hat{p}}{\partial X_2}\right).
\end{array}
\end{equation}
Using the definitions of  the $C^{\alpha}(\overline{\Omega})$ and $H^{1+\alpha,
\alpha}(\overline{\Omega}_T)$  norms, we  deduce  the
continuity of   $G$
$$\Big \| \frac{\partial\hat{p}}{\partial X_i} \Big \| _{C^{\alpha}(\overline{\Omega})}\leq \|\hat{p}\|_{H^{1+\alpha,
\alpha}(\overline{\Omega}_T)}.$$

\item $Q$ is the operator that  applies   each $ W \in
(C^{\alpha}(\overline{\Omega}_T))^2$ in $Q(W)\in
C^{\alpha}(\overline{\Omega}_T)$,  that   is the solution of the
semilinear parabolic problem \eqref{28}, given in Lemma
\ref{lem-par}. $Q$ is a continuous operator when restricted to $G
\circ D \circ I \circ E(B_M(0))$. In fact, if  $W^n :=  (W_1^n,W_2^n)$ is a sequence  that convergences to $W := (W_1 ,W_2 )$ in $(C^{\alpha}(\overline{\Omega}_T))^2$, when  $n\rightarrow \infty$,  then  also the corresponding  images by the operator $Q$ converge, that is,
$$C^n :=   Q\big((W_1^n,W_2^n)\big) \longrightarrow C:= Q\big ((W_1 ,W_2 )\big ) \quad  \mbox{in} \; C^{\alpha}(\overline{\Omega}_T),
 \;  \mbox{when} \; n\rightarrow \infty. $$

This statement can be proved  using {\it a priori} estimates for linear
parabolic problems. In effect, defining   $U^n=C^n-C$  and
$$f^n(X,t)
=(\widetilde{ \beta m}-\widetilde{\gamma
m})U^n(1-(C^n+C))+\widetilde{A_{ij}m} (W^n-W) \frac{\partial
C^n}{\partial X_i}$$
it can be derived (easily) that $U^n$  is the solution of
\begin{equation}\label{continuity}
\left \{
\begin{array}{lcl}
\displaystyle\frac{\partial U^n}{\partial t}  -  \widetilde{A_{ij}m} W
\frac{\partial U^n}{\partial X_i} -
D\widetilde{A_{ij}m}\frac{\partial^2 U^n}{\partial X_i
\partial X_j}  =   f^n(X,t) &  \mbox{in}  &   \Omega \times ]0,T]\\
  \\
 U^n  =  0  & \mbox{on}  & \partial \Omega \times ]0,T]\\
  \\
U^n(X,0)  =  0 &   \mbox{in} &  \Omega \times \{0\}.
\end{array}
\right.
\end{equation}

As $f^n$ verifies the hypotheses of \cite[Theorems 6 and 7, p. 65]{fri2}, then the next estimate for $
\|U^n\|_{C^{\alpha}(\overline{\Omega}_T)}$ is valid
$$ \|U^n\|_{C^{\alpha}(\overline{\Omega}_T)}  \leq \Big  \|(\widetilde{ \beta m}-\widetilde{\gamma
m})U^n(1-(C^n+C))+\widetilde{A_{ij}m} (W^n-W) \frac{\partial
C^n}{\partial X_i} \Big \|_{C^{\alpha}(\overline{\Omega}_T)}.$$

Moreover since $\| C\|_{C^{\alpha}(\overline{\Omega}_T)}$, $\|
C^n\|_{C^{\alpha}(\overline{\Omega}_T)}$ and $\| \frac{\partial
C^n}{\partial X_i}\|_{C^{\alpha}(\overline{\Omega}_T)}$ are bounded
we can conclude that for some  constant $K_1$,   the following relation is valid
$$ \|C^n-C\|_{C^{\alpha}(\overline{\Omega}_T)}  =  \|U^n\|_{C^{\alpha}(\overline{\Omega}_T)}   \leq  K_1 \left (|(\widetilde{ \beta m}-\widetilde{\gamma
m}) |\|C^n-C\|_{C^{\alpha}(\overline{\Omega}_T)}+
\|(W^n-W)\|_{C^{\alpha}(\overline{\Omega}_T)}\right ). $$
Supposing that
$|(\widetilde{ \beta m}-\widetilde{\gamma m})| $  is small enough we
obtain
$$\|C^n-C\|_{C^{\alpha}(\overline{\Omega}_T)}  \leq
K_2\|(W^n-W)\|_{C^{\alpha}(\overline{\Omega}_T)}$$
where $K_2$ is another constant. This proves  the sequential continuity
and consequently  the continuity of $Q$.

\end{itemize}
Finally, for $ M>  0$ and with   the conditions of
\cite[Theorem 8, p.204]{fri2}, we must show that $\mathcal {F}(B_M (0))\subset B_M (0)$,
in order to apply Schauder's fixed point theorem. As $ G \circ D \circ I $ is  a linear continuous operator,  then
\begin{equation}\label{inclusion}
\|G \circ D \circ I  \big ( E(C)\big )\|_{(C^{\alpha}(\overline{\Omega}_T))^2}  \leq
 K_ 3 \|E(C)\|_{C^{\frac{\alpha}{2}}([0,T];C^{1,
\alpha}(\overline{\Omega}))}
\end{equation}
with $K_3$ a continuity constant. Taking into account  Remark \ref{rem2},  the
 second term  of
\eqref{inclusion} can be chosen  arbitrarily small.  In this way,  it is possible to obtain the right
upper bound for the coefficients  of the convective term of the
parabolic operator, see \eqref{fried_conition}, and consequentely Lemma \ref{lem-par} ensures that $\mathcal {F}(C) \in
B_M (0)$. Therefore  the solution of  the homogenization model
\eqref{eq:homog} is a fixed point of $\mathcal {F}$. $\Box$

\section{Approximate solution of the homogenization model}\label{sec:approx}

In this section,  a numerical approach  to solve the
homogenized system \eqref{eq:homog} is briefly described.
The space-time discretization is based
on continuous finite elements for the space variable, and finite differences (using an
implicit backward Euler scheme) for the time variable.

The variational formulation corresponding to the  elliptic equation in
\eqref{eq:homog} is

\begin{equation} \label{eq:varhomog1}
\int_{\Omega} \widetilde{A_{ij}m} \frac{\partial p^0}{\partial X_j}
\frac{\partial v}{\partial X_i} dX_1 dX_2 = \int_{\Omega}[
 ( \widetilde{\beta m} -\widetilde{\gamma  m}) \,
 C^{0} + \widetilde{\gamma m}  ] \, v  dX_1 dX_2, \,\,\,\,\,\,
\forall v \in H_0^1(\Omega).
\end{equation}

  Then, denoting by  $\{\varphi_k\} $ the global finite element shape functions,  and
 considering a partition in the time interval $[0,T]$, such that
 $[0,T]=\cup_{n=0}^{N-1} [t_n,t_{n+1}]$ with  time step size $ \Delta t =t_{n+1}-t_n$, \eqref{eq:varhomog1}  leads  to the
 following approximation of the homogenized elliptic problem
\begin{equation}\label{eqhp_discr}
    \widetilde{G}_m  {\bf p^n} = (\widetilde{\beta m} -\widetilde{\gamma  m} )\, M \, {\bf C^n} \,+
 \, {\bf \widetilde { \boldsymbol{ \gamma  m }  }} .
\end{equation}
Here ${\bf p^n} = (p_k^n)$ and ${\bf C^n} = (C_k^n)$ are vectors, corresponding to the space-time discretizations of $p^0$ and $C^0$, respectively
(${\bf C^n}$ depends explicitly and ${\bf p^n}$ implicitly  on the time variable $t$). The upper index ($\bf  n$ for vector and $n$ for  vector component) corresponds to the approximation at time $t_n$  and the lower index  $k$  to the
finite element node number considered. Finally, $M $ and $\widetilde{G} $  are, respectively,   the finite element mass matrix and modified stiffness matrix  and  ${\bf \widetilde {  \boldsymbol{\gamma m}   }} =  (\widetilde {  \gamma m   }_k)$ is  a vector, defined by
\begin{equation}\label{G1m}
\begin{array}{lcl}
 {\bf \widetilde { \boldsymbol{ \gamma  m}   }} =  (\widetilde {  \gamma m   }_k) & \quad &
\widetilde {  \gamma m   }_k  :=
\displaystyle  \int_{\Omega} \widetilde {  \gamma m   } \, \varphi_k  dX_1 dX_2,
\\ 
\\
M=(M_{k,l}) &\quad & M_{k,l} :=
\displaystyle  \int_{\Omega} \varphi_k \varphi_l dX_1 dX_2,\\
\\
  \widetilde{G}_m  = (\widetilde{G}_{k,l}) &\quad &  \widetilde{G}_{k,l}:= \displaystyle  \int_{\Omega}
\nabla\varphi_k \cdot (\widetilde{Am}\nabla\varphi_l) dX_1 dX_2,
\end{array}
 \end{equation}
  with $\widetilde{Am}=\left( \widetilde{A_{ij}m} \right)_{i,j=1,2}$.

We can use similar arguments for approximating  the parabolic equation in
\eqref{eq:homog}.  Its variational formulation   is
\begin{equation}\label{eq:varhomog2}
\begin{array}{l}
  \displaystyle  \int_\Omega  \left ( \frac{\partial C^0}{\partial t} - \widetilde{A_{ij}m}\frac{\partial C^{0}}{\partial X_i} \frac{\partial
p^{0}}{\partial X_j}  \right )\, v \, dX_1 dX_2 =\\
\\
   \displaystyle  \int_\Omega   \left ( - D\widetilde{A_{ij}m}\frac{\partial  C^{0}}{ \partial X_j}  \frac{\partial  v}{ \partial X_i}+
 (\widetilde{\beta m}-\widetilde{\gamma m}) C^0(1-C^0) \, v \right )\,  dX_1 dX_2, \quad \forall v \in H_0^1(\Omega).
 \end{array}
\end{equation}
It involves, in particular,  a bilinear term in $C^0$ and $p^0$ and a non-linear term in $C^0$. For the bilinear term we remark that at time $t_n$
\begin{equation}\label{Dtilde_1step}
\begin{array}{lcl}
 \displaystyle \int_{\Omega}\widetilde{A_{ij}m} \frac{\partial C^0}{\partial X_i}
\frac{\partial p^0}{\partial X_j} v \, dX_1 dX_2
& = &  \displaystyle  \int_{\Omega}
  \frac{\partial C^0}{\partial X_i} (\widetilde{Am}\nabla p^0)_i v \, dX_1 dX_2 \\
\\
 & \simeq & \displaystyle  \sum_{s,l} p_s^n  C_l^n \int_{\Omega} (\nabla \varphi_l\cdot \widetilde{Am} \nabla \varphi_s) v \, dX_1 dX_2,
 \end{array}
\end{equation}
and then, using $v=\varphi_k$ we have
$$\int_{\Omega }\widetilde{A_{ij}m} \frac{\partial C^0}{\partial X_i}
\frac{\partial p^0}{\partial X_j} \varphi_k \, dX_1 dX_2  \simeq
(\widetilde{D}_{m,\nabla  p }  {\bf C^n})_k$$ where
\begin{equation}\label{Dgradp0}
(\widetilde{D}_{m,\nabla p})_{ k,l}:=\displaystyle \sum_{s}
p_s^n \int_{\Omega} \nabla\varphi_l\cdot
\widetilde{Am}\nabla\varphi_s \,  \varphi_k \, dX_1 dX_2.
\end{equation}
Using a backward Euler discretization in time with   time step size $\Delta
t$, the variational formulation  \eqref{eq:varhomog2}  is   approximated by the equation
\begin{equation}\label{int}
M  \frac{{\bf C^{n+1}} - {\bf C^{n}}}{\Delta t} - \widetilde{D}_{m, \nabla p}  {\bf C^{n+1}}=
-\widetilde{G}_{m,D} {\bf C^{n+1}} + \Big [ M_{\widetilde{\beta m}
-\widetilde{\gamma m}} ({\bf 1}-  {\bf C^{n }} ) \Big ] {\bf C^{n+1}}
\end{equation}
where the last term in the right hand side corresponds to the approximation of the  non-linear term  $ (1-C^0) C^0$
 (by  using  $( { \bf 1} - { \bf C^{n}})$  and ${ \bf C^{n+1 }}$ for approximating $ (1-C^0) $ and $C^0$, at time $t_n$ and $t_{n+1}$, respectively). In \eqref{int}
\begin{equation}\label{G1m_0}
\begin{array}{lcl}
M_{\widetilde{\beta m}
-\widetilde{\gamma m}}({\bf 1}-  {\bf C^{n }})=(M_{k,l}) &\quad & M_{k,l} :=
\displaystyle  \int_{\Omega} ({\widetilde{\beta m}
-\widetilde{\gamma m}} )\, (1-C_s^n) \varphi_s \, \varphi_k \, \varphi_l \, dX_1 dX_2,\\
\\
  \widetilde{G}_{ m,D} = (\widetilde{G}_{k,l}) &\quad &  \widetilde{G}_{k,l}:= \displaystyle  \int_{\Omega}
 D \nabla\varphi_k \cdot (\widetilde{Am}\nabla\varphi_l) dX_1 dX_2,
\end{array}
 \end{equation}
are other modified mass and stiffness finite element matrices.
 With a such approximation the parabolic equation  can be solved by
 computing $\bf C^{n+1}$ at each time step $t_{n+1}$, by solving  the linear system
\begin{equation}\label{eqhC_discr}
\Big (M  + \Delta t \big (-\widetilde{D}_{m,\nabla p} + \widetilde{G}_{ m,D} -
 M_{\widetilde{\beta m}-\widetilde{\gamma m}}(1-{\bf C^{n}}) \big ) \Big ){\bf C^{n+1}}=M  {\bf C^{n}}.
\end{equation}

Summarizing, the algorithm implemented in the time interval $[0,T]$,
with time step size $\Delta t= t_{n+1}-t_n$, and the initial condition ${\bf
C^{0}} $ (for the density of abnormal cells) is listed below.

\begin{itemize}
\item  Compute ${\bf p^{0}}$ the approximation of the pressure at  time $t=0$,
  by using \eqref{eqhp_discr}, {\it i.e.}
$$
\widetilde{G}_m  {\bf p^0} = 
( \widetilde{\beta m} -\widetilde{\gamma  m})\, M \, {\bf C^0}  \,  +\,  
  \widetilde{ \bf{\boldsymbol{\gamma m }}}
$$
and set $n=0$.
\item    While   $t_n < T$
\begin{enumerate}
\item   Build the matrix $\widetilde{D}_{m,\nabla p}$ defined in \eqref{Dgradp0},  for  ${\bf p^n}$, {\it i.e}
$$
(\widetilde{D}_{m,\nabla p})_{k,l}:=\displaystyle \sum_{s}
p_s^n \int_{\Omega} \nabla\varphi_l\cdot
\widetilde{Am}\nabla\varphi_s \,  \varphi_k \, dX_1 dX_2.
$$
\item Compute $\bf C^{n+1}$  using the backward Euler equation \eqref{eqhC_discr} with  $\widetilde{D}_{m,\nabla p}$
 computed in the previous step 1,   {\it i.e}
 $$
\Big (M  + \Delta t \big (-\widetilde{D}_{m,\nabla p}  + \widetilde{G}_{ m,D} -
 M_{\widetilde{\beta m}-\widetilde{\gamma m}}(1-{\bf C^{n}})\big ) \Big ){\bf C^{n+1}}=M  {\bf C^{n}}.
$$
\item Compute $\bf p^{n+1}$ using \eqref{eqhp_discr}, {\it i.e.}
$$
\widetilde{G}_m  {\bf p^{n+1}} = ( \widetilde{\beta m} -\widetilde{\gamma  m})\, M \, {\bf C^{n+1}}  \,  +\,  {\bf \widetilde{\boldsymbol{\gamma  m}}}  .
$$
 \item Increment $n$, that is $n=n+1$, and go to step 1.

\end{enumerate}
\item   Stop when $t_n=T$.

\end{itemize}

\section{Numerical simulations}\label{sec:simu}
In this section we present the numerical simulations for the homogenized model \eqref{eq:homog},
obtained by implementing the algorithm described in the previous Section \ref{sec:approx}.
Moreover we compare these results with those obtained by solving numerically the heterogeneous problem \eqref{11}.
 
    The  crypt  geometry is  related to the real dimensions  of colonic crypts:
 the domain, representing a piece of the colon,  is
   $\Omega=[-1,1]\times[-1,1]$ (a rectangle as well could have been considered,  as shown in  Figure \ref{fig:crypts}) and  the values for the radius of the shrunk crypt bottom, the radius of the top crypt orifice and the crypt height, in \eqref{1}-\eqref{2}, are  $r=a/4$, $R=a/2$ and $L=14a$, respectively,
 where $a=\frac{\sqrt{2}}{3^{\frac{3}{4}}}$ is the edge of the reference hexagon $P$ of area 1. We have used the relation $L/ 2R = 14$,  between the dimensions of the  height $L$ and the diameter $2R$, according to the dimensions reported for human colonic crypts in \cite{guebel2008computer, Halm} ($433\, \mu m$ for the crypt height, $16\, \mu m$ for the
    crypt top orifice excluding the epithelium and $15.1 \mu m$ for the epithelium cell depth).
The time unit measure used is the hour. 
 Triangular finite  elements and 3-point Gaussian quadrature rule are used for the finite element discretization used to get $m$ that satisfies     \eqref{22} in $P$. Bilinear finite elements  and 4-point Gaussian quadrature  rule are used for the finite element discretization of the systems \eqref{11} and \eqref{eq:homog}.

Concerning the values for the parameters, we define the  proliferative coefficients $\gamma $ and $\beta $, introduced in
\eqref{eq:modelS2_S3}, by
\begin{equation}\label{values}
\begin{array}{lcl}
\gamma (x_3) & = & \left\{ \begin{array}{cl} \tau_{\gamma}(x_3-2/3 L)^2 &  {\rm if } \quad x_3\leq \frac{2}{3}L\\
[.25em]
              0 & {\rm elsewhere }\end{array}\right. \\ \\
\beta (x_3) & = & \left\{ \begin{array}{cl} \tau_{\beta_1}(x_3-2/3
L)^2 + \beta_2 &
    {\rm if } \quad x_3\leq \frac{2}{3}L \\
    [.25em]
              \beta_2 & {\rm elsewhere }\end{array}\right.
\end{array}
\end{equation}
with $\tau_{\gamma }= \tau_{\beta_1}=0.01$, $\beta_2=0.1$.
   This
choice reassures the known qualitative behavior of colonic cells
 to have a high proliferation at the bottom of the crypt (where $x_3=0$) that decreases going upwards
to the orifice (where $x_3=L$) along the crypt axis,
see \cite{figu6}. Moreover, as assumed by several authors  (see for instance \cite{DL2001,JEBMC}) we suppose that there is
no proliferation in the last third of the crypt (starting from the
bottom),  that is for $x_3 \geq \frac{2}{3}L$. We use the following constant
diffusion coefficients $D_N=D_C= 0.1$ in \eqref{D1D2} that
guarantee
  \eqref{constantDs} is satisfied.  

We assume that  at an initial time (t = 0)   there exists already a non-zero density of abnormal cells, located in a region $\Omega$,  containing several  crypts. This initial   density of abnormal cells and their
  location in $\Omega$ is represented by the function
   $C(X_1,X_2,0)=0.9 e^{-\frac{X_1^2+X_2^2}{0.1}}$.
  Obviously the simulations can be done for an arbitrary initial density of abnormal cells, it is just a matter of redefining this function  $C(X_1,X_2,0)$ that represents the location and number of abnormal cells at the initial time (in particular, the abnormal cells can be located only in  a single crypt).

\begin{figure}[h!]
 \centering
 \includegraphics[scale = 0.3]{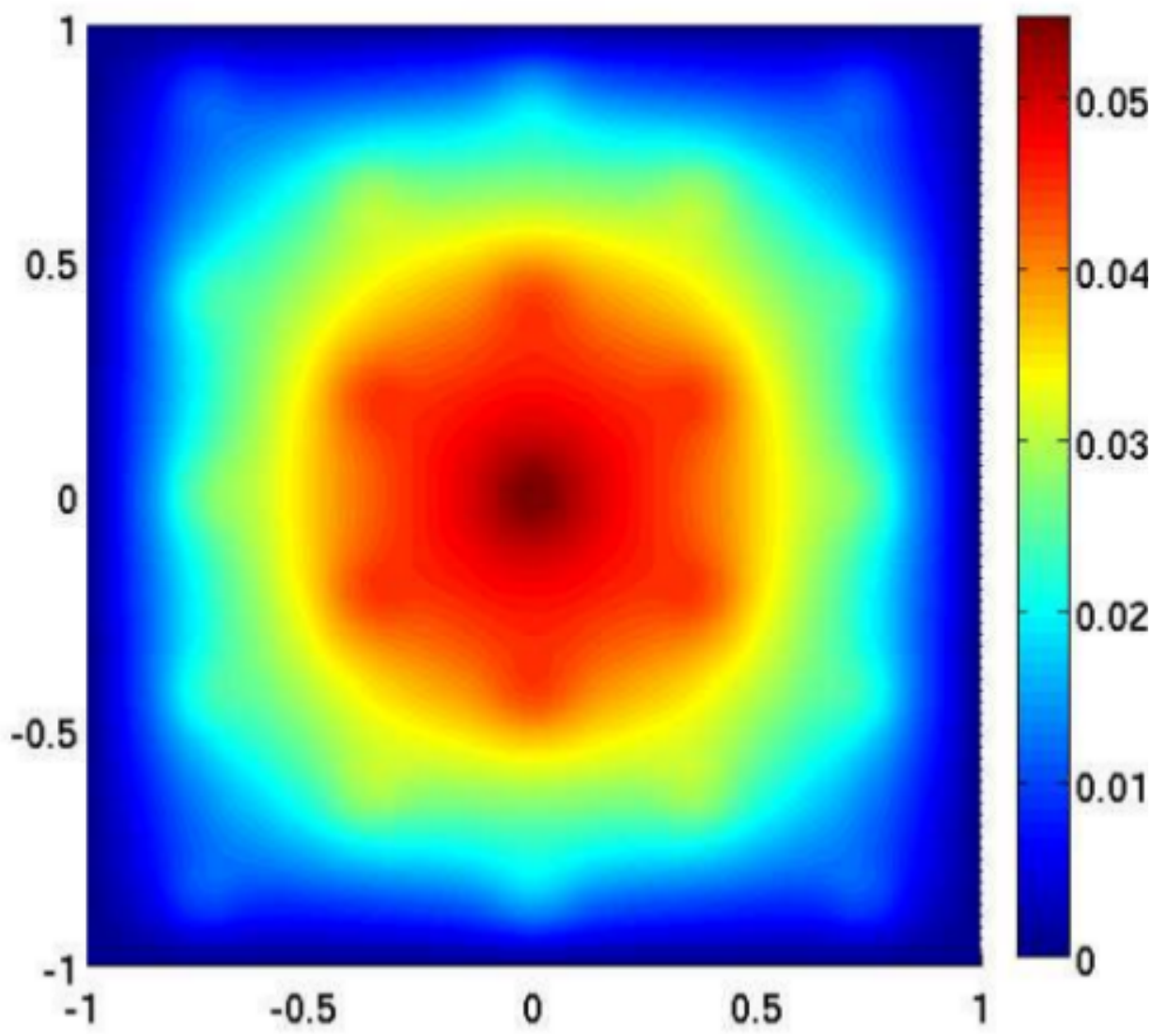}  \hspace{1cm}
  \includegraphics[scale = 0.3]{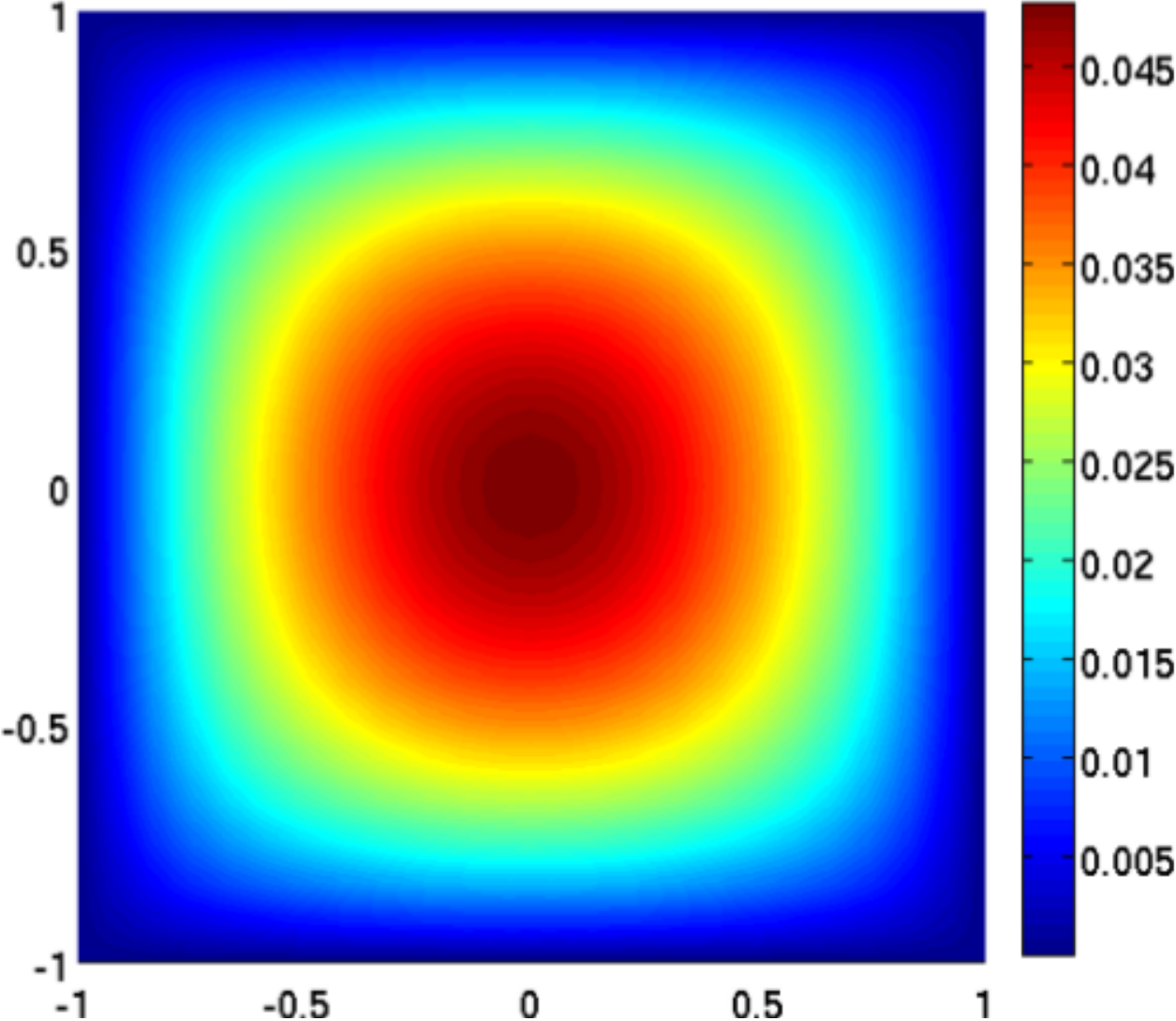}  \hspace{1cm}
\includegraphics[scale = 0.3]{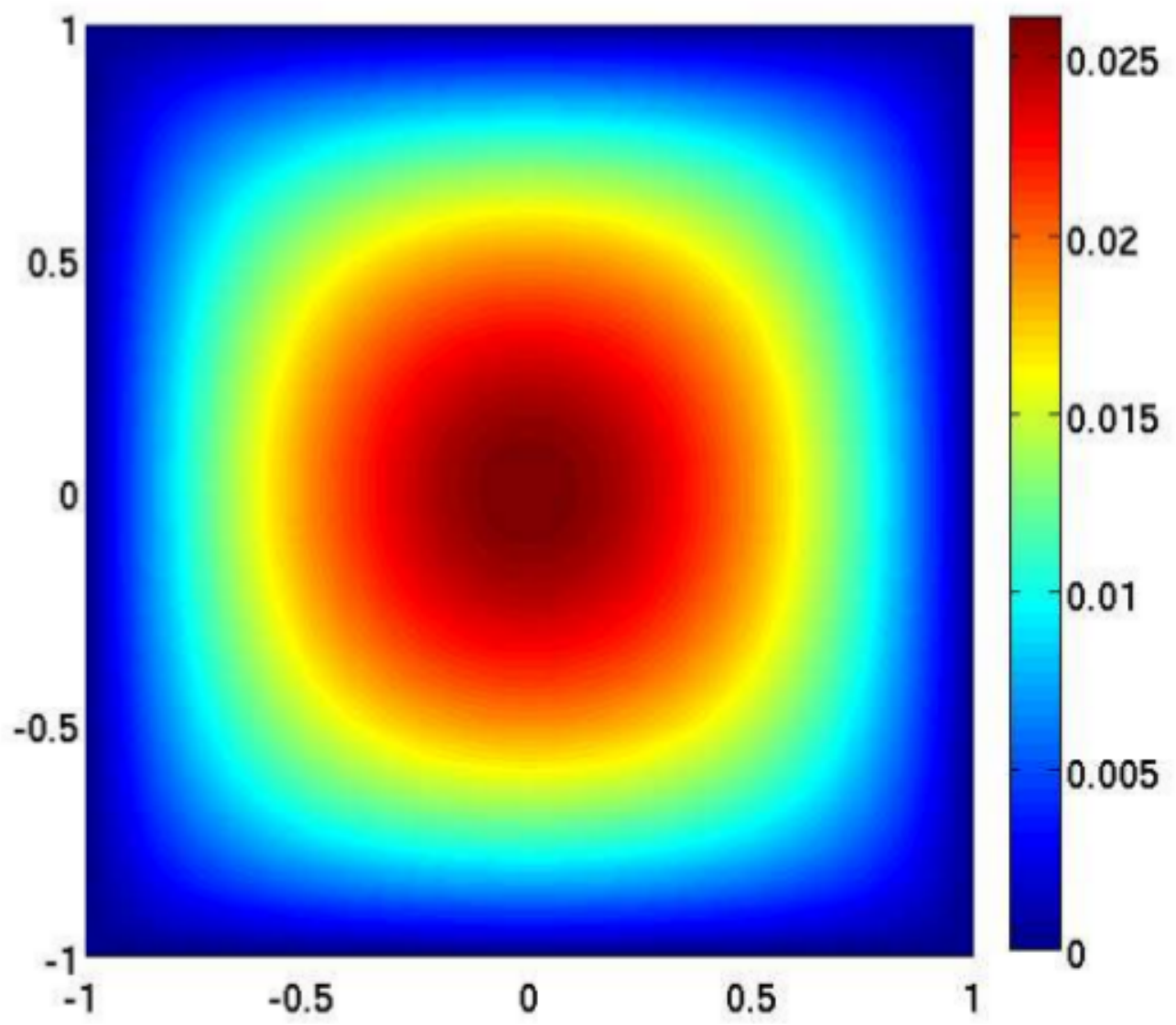} \\
[.75em]
\includegraphics[scale = 0.3]{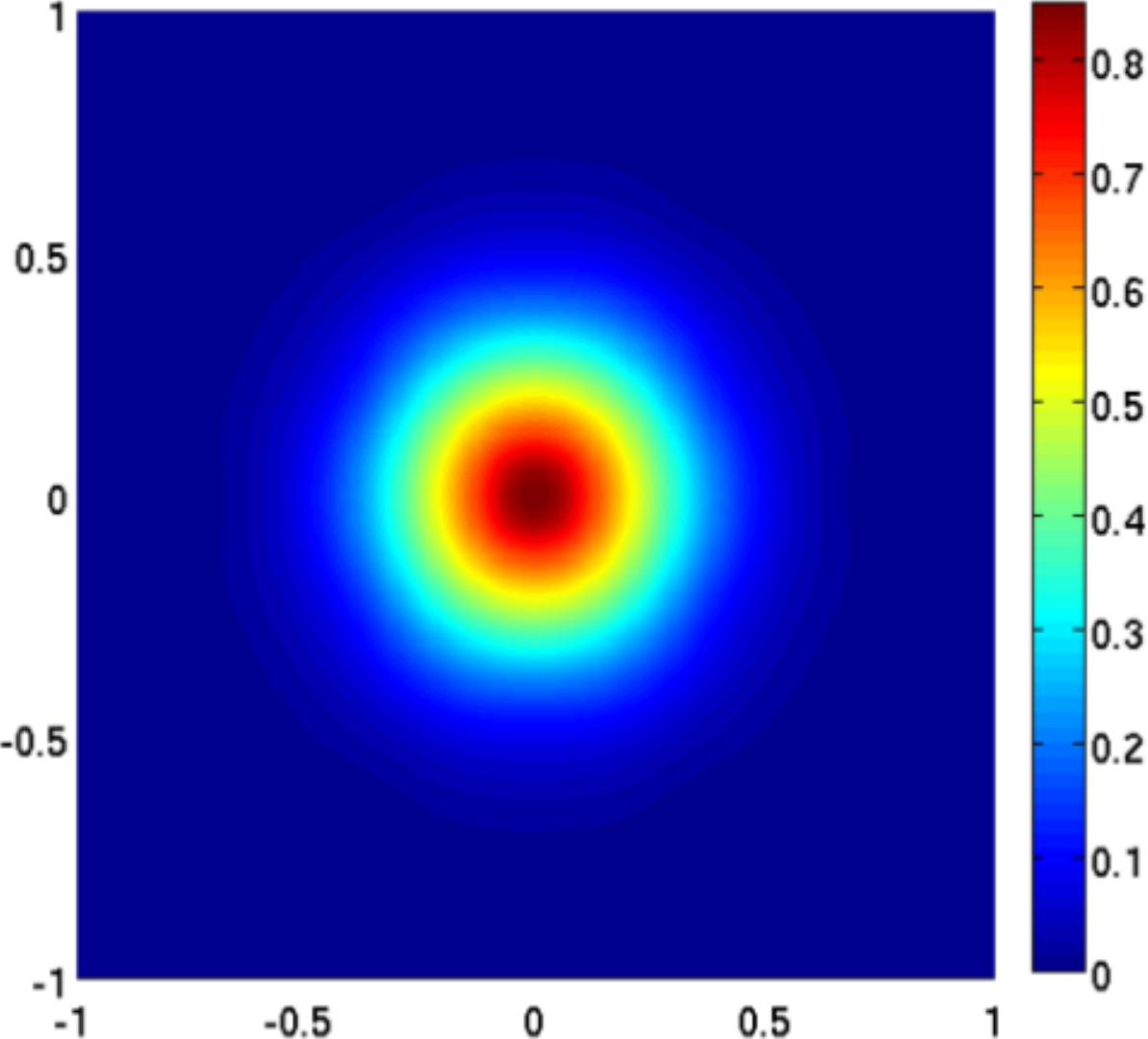}  \hspace{1cm}
  \includegraphics[scale = 0.3]{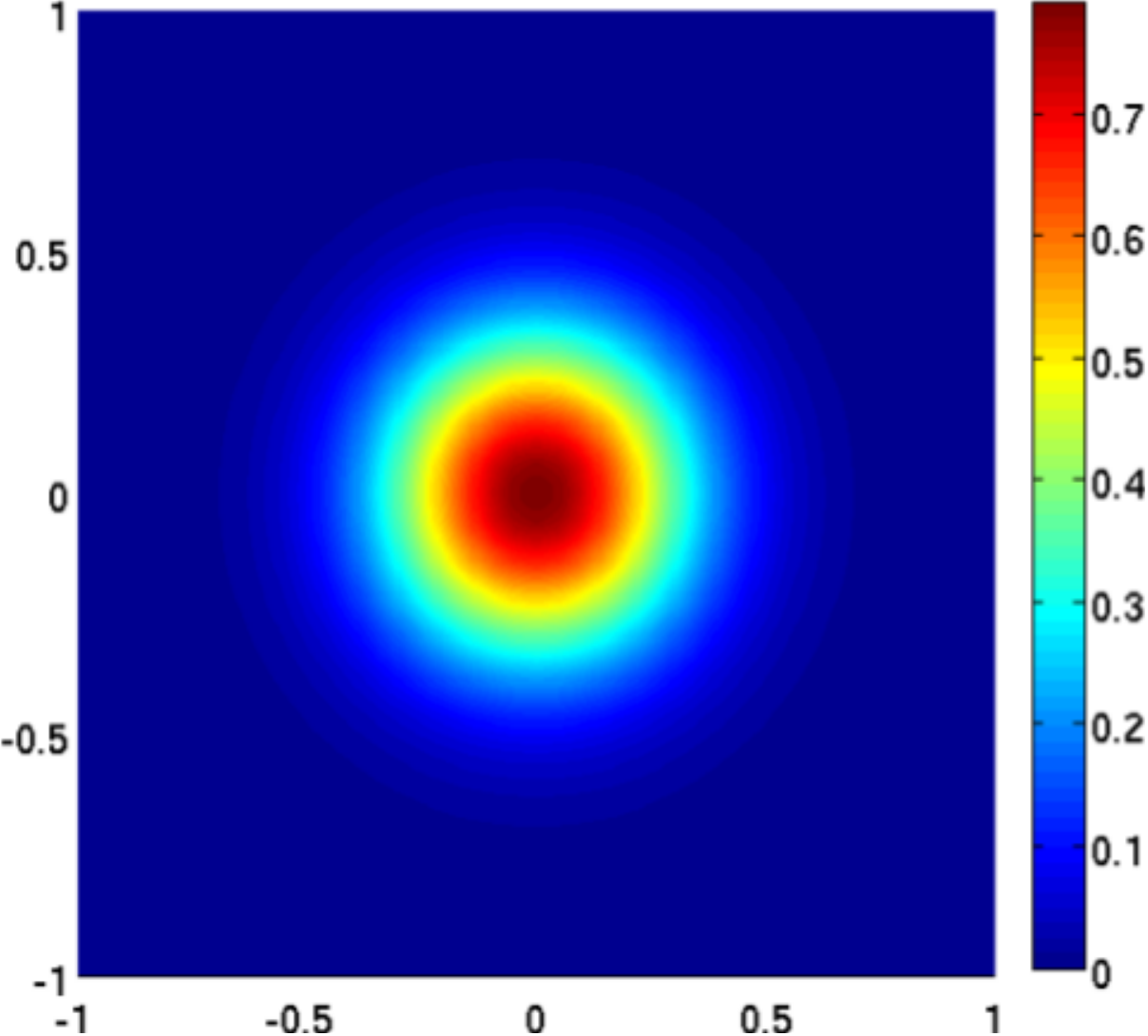} \hspace{1cm}
  \includegraphics[scale = 0.3]{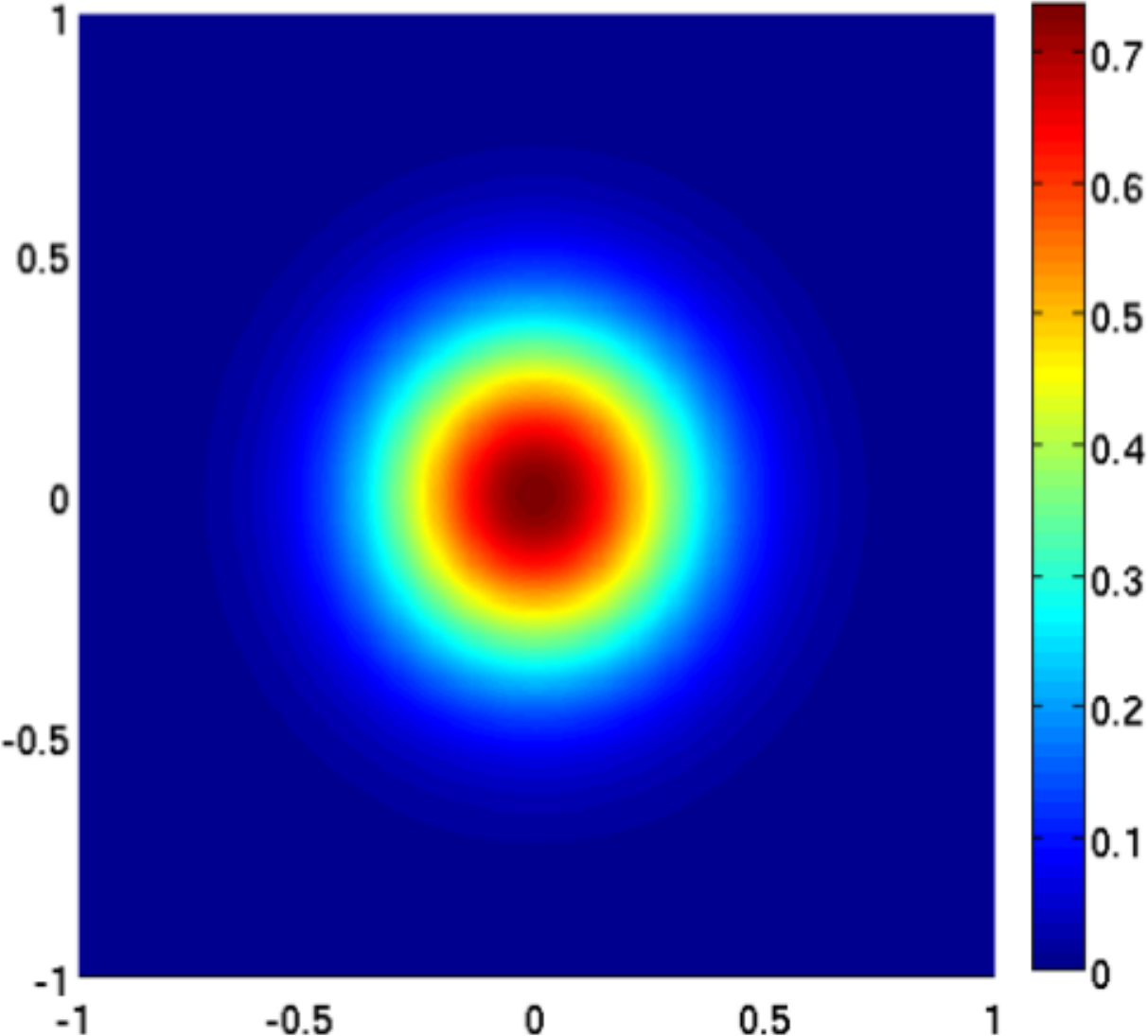}
\caption{ Left and middle columns: Solutions of the heterogenous periodic model \eqref{11} for two different values of $\varepsilon$
 -- pressure $p^{\varepsilon}$ (top)  and cell density $C^{\varepsilon}$ (bottom) at time $t=0.05$. In the left column
 $\varepsilon=0.4$ 
 and  in the middle column
 $\varepsilon=0.05$. 
 Right column: Solution of the homogenized model  \eqref{eq:homog} -- pressure  $p^0$ (top) and cell density $C^0$ (bottom) at time  $t=0.05$.}
\label{figp}
\end{figure}

The fine scale solution $\big ((C^{\varepsilon}(X,t), p^{\varepsilon}(X,t) \big)$, of the heterogeneous periodic model \eqref{11}, is
expected to converge to $\big ((C^{0}(X,t),
p^{0}(X,t) \big)$, the solution of the homogenization  problem \eqref{eq:homog}, when $\varepsilon$   goes to zero.
 The numerical results obtained,  that are presented in the  Tables \ref{tabprL2}-\ref{tabCrLinf} and Figure \ref{figp}, show this convergence. For each  fixed $\varepsilon$,  in  the experiments performed for the heterogenous periodic model \eqref{11}  the number of 
 subdomains
 periodically distributed  in the domain
  $\Omega=[-1,1]^2$ is equal to $\frac{4}{\varepsilon^2}$.

  \begin{table}[h!]
   
     \begin{tabular*}{\hsize}{@{\extracolsep{\fill}}|l|l|l|l|l|@{}}
  \hline
    $\|p^{\varepsilon} - p^0\|_{2}/\|p^0\|_{2}$ & $\varepsilon=8e-01$ & $\varepsilon=4e-01$ & $\varepsilon=5e-02$ & $\varepsilon=3.2e-02$ \\
  \hline
    t=1e-02    & 1.1230 & 9.5942e-01 & 8.4966e-01 & 7.2541e-01   \\
  \hline
    t=3e-02    & 1.2218 & 9.5817e-01 & 8.4836e-01 & 7.2412e-01  \\
  \hline
    t=5e-02    & 1.2119 & 9.5668e-01 & 8.4692e-01 & 7.2268e-01 \\
  \hline
  \end{tabular*}\caption{Relative error for pressure in $L^{2}$ norm at times
$0.01,0.03,0.05$ with space step size $h=5e-03$ and time step size
$\Delta t=5e-03$. }\label{tabprL2}

  \end{table}
\begin{table}[h!]
  
  \begin{tabular*}{\hsize}{@{\extracolsep{\fill}}|l|l|l|l|l|@{}}
 \hline
    $\|p^{\varepsilon} - p^0\|_{\infty}/\|p^0\|_{\infty}$ & $\varepsilon=8e-01$ & $\varepsilon=4e-01$  & $\varepsilon=5e-02$ & $\varepsilon=3.2e-02$  \\
 \hline
   t=1e-02  & 1.6362 & 1.0903 & 8.4294e-01 & 7.2272e-01  \\
 \hline
   t=3e-02  & 1.6431 & 1.0928 & 8.4455e-01 & 7.2385e-01 \\
 \hline
   t=5e-02  & 1.6472 & 1.0942 & 8.4546e-01 & 7.2434e-01 \\
 \hline
 \end{tabular*} \caption{Relative error for pressure in $L^{\infty}$ norm at times $0.01,0.03,0.05$ with space step size $h=5e-03$ and time step size
          $\Delta t=5e-03$. }
\label{tabprLinf}

 \end{table}

 Tables \ref{tabprL2}-\ref{tabCrLinf}   list the relative
errors
 between the fine-scale solution $ (C^{\varepsilon},p^{\varepsilon})$ and the homogenized solution $ (C^{0},p^{0})$.
  In particular the relative errors for the pressure  and cell density  in the  $L^2$ and
$L^{\infty}$ norms (denoted by $\|.\|_2$ and $\|.\|_\infty$) are listed, respectively,  in Tables \ref{tabprL2},
\ref{tabprLinf}  and
Tables \ref{tabCrL2},  \ref{tabCrLinf}. We can see that halving the
parameter $\varepsilon$ the relative errors for the pressure and cell density
decrease. This shows, numerically, that $(C^{\varepsilon},p^{\varepsilon})$ becomes closer
to  $ (C^{0},p^{0})$ when $\varepsilon$ decreases.

 \begin{table}[h!]
 
 \begin{tabular*}{\hsize}{@{\extracolsep{\fill}}|l|l|l|l|l|@{}}
 \hline
  $\|C^{\varepsilon} - C^0\|_{2}/\|C^0\|_2$ & $\varepsilon=8e-01$ & $\varepsilon=4e-01$  & $\varepsilon=5e-02$ & $\varepsilon=3.2e-02$  \\
 \hline
   t=1e-02   & 2.0961e-02 &  1.5200e-02   & 1.2808e-02   & 1.2595e-02 \\
 \hline
   t=3e-02   & 5.8654e-02 &  4.1922e-02   & 3.5728e-02   & 3.4701e-02 \\
 \hline
   t=5e-02   & 9.1470e-02 &  6.4395e-02   & 5.4946e-02   & 5.3155e-02 \\
 \hline 
 \end{tabular*}
 \caption{Relative error for cell density in $L^{2}$ norm at times $0.01,0.03,0.05$ with space step size $h=5e-03$ and time step size $\Delta t=5e-03$.}
\label{tabCrL2}

 \end{table}

\begin{table}[h!]
 
    \begin{tabular*}{\hsize}{@{\extracolsep{\fill}}|l|l|l|l|l|@{}}
 \hline
  $\|C^{\varepsilon} - C^0\|_{\infty}/\|C^0\|_{\infty}$ & $\varepsilon=8e-01$ & $\varepsilon=4e-01$  & $\varepsilon=5e-02$ & $\varepsilon=3.2e-02$  \\
\hline
   t=1e-02   & 3.6303e-02 &  3.6179e-02  &  2.0593e-02  & 1.9457e-02 \\
 \hline
   t=3e-02   & 1.0453e-01 &  1.0202e-02  &  5.1916e-02  & 5.0650e-02 \\
 \hline
   t=5e-02   & 1.6721e-01 &  1.5773e-01  &  7.7716e-02    & 7.5956e-02  \\
 \hline
 \end{tabular*}
 \caption{Relative error for cell density in $L^{\infty}$ norm at times $0.01,0.03,0.05$ with space step size $h=5e-03$ and time step size $\Delta t=5e-03$. }
\label{tabCrLinf}
 \end{table}

Note that it is not possible
to use smaller values for $\varepsilon$ than those presented, since
 this would require the use of a small  spatial step size $h$, in the finite element mesh, and this leads  to an unacceptable high memory cost. On the other hand, since the  hexagonal distributed domain $\varepsilon P$ in $\Omega$
 has two concentric circles, the  finite elements used must satisfy
 a condition stronger than just $h < \varepsilon$ to have high accuracy in the numerical approximation
 of the solution of \eqref{11}.
In our case, since we use  quadrilateral finite elements  with a regular discretization  for $\varepsilon P$,
 there are  finite elements inside each  circle of $\varepsilon P$  only if   only if   $h\leq \frac{\varepsilon a}{4}$ 
(in $\varepsilon P$ the inner and outer circles have radius equal to $\frac{\varepsilon a}{4}$ and $\frac{\varepsilon a}{2}$, respectively).
 In particular, the value of $\varepsilon$ used in the last columns  of  Tables \ref{tabprL2}-\ref{tabCrLinf} is $\frac{4h}{a}\approx 3.2e-02$, 
that is the minimum value  for  $\varepsilon$   in order to have finite elements of size $h$ inside each region of $\varepsilon P$.

 Figure \ref{figp} shows a visible homogenization of the pressure when $\varepsilon$ decreases. We observe in fact that
$p^{\varepsilon}$ has oscillations for high $\varepsilon$, that disappear when $\varepsilon$ decreases.
The density $C^\varepsilon$  shows instead a slight
visible convergence to $C^0$. A convergence can only be noticed by
the fact that the maximum value $C^\varepsilon$
tends to the maximum value of $C^0$ when $\varepsilon$ decreases as seen in the toolbars of the plots in Figure \ref{figp}.

Finally we have compared  the computational
time required to obtain the numerical approximation of the periodic
multiscale problem \eqref{11} and of the homogenization problem
\eqref{eq:homog}. The former takes approximately   1800 seconds  to
execute a single iteration in time (to obtain $ (C^{\varepsilon},p^{\varepsilon})$)
 with a spatial step size $5e-03$ in $\Omega$, and this  is more than seven times larger than the time spent by the homogenized problem,
that needs in fact just 250 seconds to obtain  $ (C^{0},p^{0})$. The main reason
for this computational time  difference, between the two problems, is due to the fact that
the coefficients of the homogenized problem \eqref{eq:homog} are
constants  with respect to the spatial variable, whereas this is not the case for  coefficients of the  periodic multiscale problem \eqref{11}, and then the latter problem needs a larger time cost.

\section{Conclusions}\label{sec:conc}

In this paper we derived an homogenization model, for simulating the ACF spread and dynamics,  with the goal of reproducing {\it in silico} at a macroscopic level
({\it i.e.} at the tissue level) what is observed in conventional colonoscopy images.
We have  considered the  PDEs model \eqref{11} for  the mathematical description of the biological phenomenon. It is an heterogenous periodic model based on the assumption that the colon is a domain, with a periodic distribution of very small heterogeneities, the crypts,  wherein most cellular activity occurs.
In this model we suppose the existence of an abnormal cell population that has the property to proliferate not only inside the crypts,  as done by normal cells, but also outside the crypts. This property is believed to be responsible for abnormal cells to invade neighbor crypts and possibly  induce crypt fission and the formation of an adenoma. Then we have applied an  homogenization technique to the  heterogenous periodic model that relies on the two-scale asymptotic expansion method.

In a previous work \cite{figu6} we solved numerically an heterogeneous periodic problem, similar to \eqref{11}, without any homogenization technique, but applying a multi-scale numerical approach. In \cite{figu6} we have observed that choosing a small size for the heterogeneousness ({\it i.e.} the crypt size $\varepsilon$) results in a high cost of memory and computations.

 The main benefit of the homogenization model herein proposed, in comparison  with the approach developed in \cite{figu6},  is that it permits to describe with a simpler model (the homogenized one) a very complex, periodic and multi-scale problem. For the homogenization model, a numerical solution can be easily computed, by means of standard discretization methods, as finite elements  for the spatial variable combined with finite differences for the time variable.

  In addition the advantage of this homogenization model is the possibility to perform experiments {\it in silico}, that represent simulations of the ACF evolution at a macroscopic, or equivalently, tissue level,  starting with an arbitrary density (small or big) of abnormal cells that can be located in an arbitrary colon region (inside a single crypt or in several crypts, adjacent or not, or outside the crypts). In effect, the location and density of the abnormal cells at the initial time is controlled (meaning represented) by the definition of the initial cell density.

 In the future we intend to build a more complete hybrid model, by incorporating   other biological effects of the ACF formation,   such as the crypt fission and the viscoelastic properties of the colon epithelium. This will permit to compare the predicted ACF dynamics with that obtained by medical images and histopathological results.

\section*{Acknowledgment}
 The authors   thank the anonymous referees for
valuable comments and suggestions that helped to improve the
manuscript.

\bibliographystyle{plain}
\bibliography{crypts_refs}

\begin{thebibliography}{10}

\bibitem{Araki1996}
K.~Araki, Y.~Furuya, M.~Kobayashi, K.~Matsuura, T.~Ogata, and H.~Isozaki.
\newblock Comparison of microvascularture between the proximal and distal human
  colon.
\newblock {\em J. Electron. Microsc.}, 45:202--206, 1996.

\bibitem{APS}
N.~J. Armstrong, K.~J. Painter, and J.~A. Sherratt.
\newblock A continuum approach to modelling cell--cell adhesion.
\newblock {\em Journal of Theoretical Biology}, 243:98--113, 2006.

\bibitem{Bakeretal}
A.~M. Baker et~al.
\newblock Quantification of crypt and stem cell evolution in the normal and
  neoplastic human colon.
\newblock {\em Cell Reports}, 8:940--947, 2014.

\bibitem{BGM2008}
R.~E. Baker, E.~A. Gaffney, and P.~K. Maini.
\newblock Partial differential equations for self-organization in cellular and
  developmental biology.
\newblock {\em Nonlinearity}, 21(11):R251--R290, 2008.

\bibitem{bakhvalov1989homogenisation}
N.S. Bakhvalov and G.~P. Panasenko.
\newblock {\em Homogenisation: averaging processes in periodic media:
  mathematical problems in the mechanics of composite materials}.
\newblock Kluwer Academic Publishers, 1989.

\bibitem{bensoussanasymptotic}
A.~Bensoussan, J.L. Lions, and G.~Papanicolaou.
\newblock {\em Asymptotic Analysis for Periodic Structures}.
\newblock North-Holland, Amsterdam, 1978.

\bibitem{BC}
M.~Bienz and H.~Clevers.
\newblock Linking colorectal cancer to {Wnt} signaling.
\newblock {\em Cell}, 103:311--320, 2000.

\bibitem{birdI}
R.~P. Bird.
\newblock Role of aberrant crypt foci in understanding the pathogenesis of
  colon cancer.
\newblock {\em Cancer Letters}, 93:55--71, 1995.

\bibitem{birdII}
R.~P. Bird and C.~K. Good.
\newblock The significance of aberrant crypt foci in understanding the
  pathogenesis of colon cancer.
\newblock {\em Toxicology Letters}, 112-113:395--402, 2000.

\bibitem{Boman2008}
B.~Boman, J.~Fields, K.~Cavanaugh, A.~Guetter, and O.~Runquist.
\newblock How dysregulated colonic crypt dynamics cause stem cell
  overpopulation and initiate colon cancer.
\newblock {\em Cancer Research}, 68(9):3304--3313, 2008.

\bibitem{carulli2014unraveling}
A.~J. Carulli, L.~C. Samuelson, and S.~Schnell.
\newblock Unraveling intestinal stem cell behavior with models of crypt
  dynamics.
\newblock {\em Integrative Biology}, 6(3):243--257, 2014.

\bibitem{cioranescu1999introduction}
D.~Cioranescu and P.~Donato.
\newblock {\em An introduction to homogenization}, volume~26.
\newblock Oxford University Press Oxford, 1999.

\bibitem{CostaSole}
J.~Costa, I.~Gonzalez-Garcia, and R.~Sole.
\newblock Spatial dynamics in cancer.
\newblock In {\em Complex Systems Science in Biomedicine, Chapter 6.2}, pages
  557--577. (Deisboeck, T.S., and Kresh, J.Y.), 2006.

\bibitem{CLLW}
V.~Cristini, X.~Li, J.~Lowengrub, and S.~M. Wise.
\newblock Nonlinear simulations of solid tumor growth using a mixture model:
  invasion and branching.
\newblock {\em Journal of Mathematical Biology}, 58:723--763, 2009.

\bibitem{CLN}
V.~Cristini, J.~Lowengrub, and Q.~Nie.
\newblock Nonlinear simulation of tumor growth.
\newblock {\em Journal of Mathematical Biology}, 46:191--224, 2003.

\bibitem{DL2001}
D.~Drasdo and M.~Loeffer.
\newblock Individual-based models to growth and folding in one-layered tissues:
  intestinal crypts and early development.
\newblock {\em Nonlinear Analysis}, 47:245--256, 2001.

\bibitem{EC}
C.~M. Edwards and S.J. Chapman.
\newblock Biomechanical modelling of colorectal crypt budding and fission.
\newblock {\em Bulletin of Mathematical Biology}, 69(6):1927--1942, 2007.

\bibitem{Eng}
B.~Engquist and P.E. Souganidis.
\newblock Asymptotic and numerical homogenization.
\newblock {\em Acta Numerica}, 17:147--190, 2008.

\bibitem{evans1998partial}
L.~C. Evans.
\newblock {\em Partial differential equations. Graduate studies in
  mathematics}, volume~19.
\newblock 1998.

\bibitem{figueiredo2013physiologic}
I.~N. Figueiredo and C.~Leal.
\newblock Physiologic parameter estimation using inverse problems.
\newblock {\em SIAM Journal on Applied Mathematics}, 73(3):1164--1182, 2013.

\bibitem{figu4}
I.~N. Figueiredo, C.~Leal, T.~Leonori, G.~Romanazzi, P.~N. Figueiredo, and M.M.
  Donato.
\newblock A coupled convection-diffusion level set model for tracking
  epithelial cells in colonic crypts.
\newblock {\em Procedia Computer Science}, 1(1):955--963, 2010.

\bibitem{figu5}
I.~N. Figueiredo, C.~Leal, G.~Romanazzi, B.~Engquist, and P.~N. Figueiredo.
\newblock A convection-diffusion-shape model for aberrant colonic crypt
  morphogenesis.
\newblock {\em Computing and Visualization in Science}, 14(4):157--166, 2011.

\bibitem{figu6}
I.~N. Figueiredo, G.~Romanazzi, C.~Leal, and B.~Engquist.
\newblock A multiscale model for aberrant crypt foci.
\newblock {\em Procedia Computer Science}, 18:1026--1035, 2013.

\bibitem{figu2010}
P.~Figueiredo and M.~Donato.
\newblock Cyclooxygenase-2 is overexpressed in aberrant crypt foci of smokers.
\newblock {\em European Journal of Gastroenterology \& Hepatology},
  22(10):1271, 2010.

\bibitem{figu2009}
P.~Figueiredo, M.~Donato, M.~Urbano, H.~Goul{\~a}o, H.~Gouveia, C.~Sofia,
  M.~Leit{\~a}o, and Diniz Freitas.
\newblock Aberrant crypt foci: endoscopic assessment and cell kinetics
  characterization.
\newblock {\em International Journal of Colorectal Disease}, 24(4):441--450,
  2009.

\bibitem{FBC}
A.~Fletcher, C.~Breward, and S.~Chapman.
\newblock Mathematical modeling of monoclonal conversion in the colonic crypt.
\newblock {\em Journal of Theoretical Biology}, 300:118--133, 2012.

\bibitem{fri1}
A.~Friedman.
\newblock On quasi-linear parabolic equations of the second order.
\newblock {\em Journal of Mathematics and Mechanics}, 7:793--808, 1958.

\bibitem{fri2}
A.~Friedman.
\newblock {\em Partial differential equations of parabolic type}.
\newblock Prentice-Hall, Englewood Cliffs, NJ., 1964.

\bibitem{frie}
A.~Friedman.
\newblock Free boundary problems arising in tumor models.
\newblock {\em Atti Accad. Naz. Lincei Cl. Sci. Fis. Mat. Natur. Rend. Lincei
  (9) Mat. Appl.}, 15(3-4):161--168, 2004.

\bibitem{greavesal}
L.~C. Greaves et~al.
\newblock Mitochondrial {DNA} mutations are established in human colonic stem
  cells, and mutated clones expand by crypt fission.
\newblock {\em Proceedings of the National Academy of Sciences of the United
  States}, 103(3):714--719, 2006.

\bibitem{Greenspan}
H.~P. Greenspan.
\newblock On the growth and stability of cell cultures and solid tumors.
\newblock {\em Journal of Theoretical Biology}, 56(2):229?--242, 1976.

\bibitem{guebel2008computer}
D.~V. Guebel and N.~V. Torres.
\newblock A computer model of oxygen dynamics in human colon mucosa:
  Implications in normal physiology and early tumor development.
\newblock {\em Journal of theoretical biology}, 250(3):389--409, 2008.

\bibitem{Halm}
P.~R. Halm and S.~T. Halm.
\newblock Secretagogue response of globet cells and columnar cells in human
  colonic crypts.
\newblock {\em Am J Physiol Cell Physiol}, 278:C212--C233, 2000.

\bibitem{HJ}
P.~R. Harper and S.~K. Jones.
\newblock Mathematical models for the early detection and treatment of
  colorectal cancer.
\newblock {\em Health Care Management Science}, 8:101--109, 2005.

\bibitem{HillMA}
M.A. Hill.
\newblock (2016) {E}mbryology gastrointestinal tract - colon histology.
  {R}etrieved {F}ebruary 29, 2016, from
  \url{https://embryology.med.unsw.edu.au/embryology/index.php/File:Colon_histology_005.jpg}.

\bibitem{hurlstoneal}
D.~Hurlstone et~al.
\newblock Rectal aberrant crypt foci identified using
  high-magnification-chromoscopic colonoscopy: biomarkers for flat and
  depressed neoplasia.
\newblock {\em American Journal of Gastroenterology}, pages 1283--1289, 2005.

\bibitem{Isele}
W.~Isele and H.~Meinzer.
\newblock Applying computer modeling to examine complex dynamics and pattern
  formation of tissue growth.
\newblock {\em Computational and Biomedical Research}, 31:476--494, 1998.

\bibitem{JEBMC}
M.~D. Johnston, C.~M. Edwards, W.~F. Bodmer, P.~K. Maini, and S.~J. Chapman.
\newblock Mathematical modeling of cell population dynamics in the colonic
  crypt and in colorectal cancer.
\newblock {\em Proceedings of the National Academy of Sciences of the United
  States}, 104(10):4008--4013, 2007.

\bibitem{Kershaw}
S.~K. Kershaw, H.~M. Byrne, D.~J. Gavaghan, and J.~M. Osborne.
\newblock {Colorectal cancer through simulation and experiment}.
\newblock {\em Systems Biology, IET}, {7}({3}):{57--73}, {2013}.

\bibitem{KingFranks}
J.~R. King and S.~J. Franks.
\newblock Mathematical analysis of some multi-dimensional tissue-growth models.
\newblock {\em European Journal of Applied Mathematics}, 15:273--295, 2004.

\bibitem{KomaWang}
Wang~L. Komarova, N.
\newblock Initiation of colorectal cancer.
\newblock {\em Cell Cycle}, 3:1558--1565, 2004.

\bibitem{lamlip}
S.~A. Lamprecht and M.~Lipkin.
\newblock Migrating colonic crypt epithelial cells: primary targets for
  transformation.
\newblock {\em Carcinogenesis}, 23(11):1777--1780, 2002.

\bibitem{LBWP}
M.~Loeffler, A.~Birke, D.~Winton, and C.~Potten.
\newblock Somatic mutation, monoclonality and stochastic models of stem cell
  organization in the intestinal crypt.
\newblock {\em Theor. Biol.}, 160:471--491, 1993.

\bibitem{MILN}
F.~Michor, Y.~Iwasa, C.~Lengauer, and M.~A. Nowak.
\newblock Dynamics of colorectal cancer.
\newblock {\em Seminars in Cancer Biology}, 15:484--494, 2005.

\bibitem{michor2004linear}
F.~Michor, Y.~Iwasa, H.~Rajagopalan, C.~Lengauer, and M.~A. Nowak.
\newblock Linear model of colon cancer initiation.
\newblock {\em Cell Cycle}, 3:358--362, 2004.

\bibitem{MFMB}
G.~R. Mirams, A.~G. Fletcher, P.~K. Maini, and H.~M. Byrne.
\newblock A theoretical investigation of the effect of proliferation and
  adhesion on monoclonal conversion in the colonic crypt.
\newblock {\em Journal of Theoretical Biology}, 312:143--156, 2012.

\bibitem{MM}
P.~J. Murray, J.~W.~A. Kang, G.~R. Mirams, S.~Y. Shin, H.~M. Byrne, P.~K.
  Maini, and K.~H. Cho.
\newblock Modelling spatially regulated $\beta$-catenin dynamics and invasion
  in intestinal crypts.
\newblock {\em Biophysical Journal}, 99, 2010.

\bibitem{MWFETM}
P.~J. Murray, A.~Walter, A.~G. Fletcher, C.~M. Edwards, M.~J Tindall, and P.~K.
  Maini.
\newblock Comparing a discrete and continuum model of the intestinal crypt.
\newblock {\em Phys. Biol.}, 8:026011, 2011.

\bibitem{nowak2002role}
M.~A. Nowak, N.~L Komarova, A.~Sengupta, P.~V. Jallepalli, I-M. Shih,
  B.~Vogelstein, and C.~Lengauer.
\newblock The role of chromosomal instability in tumor initiation.
\newblock {\em Proceedings of the National Academy of Sciences},
  99(25):16226--16231, 2002.

\bibitem{Osborneetal}
J.~Osborne, A.~Walter, S.~Kershaw, G.~Mirams, A.~Fletcher, P.~Pathmanathan,
  D.~Gavaghan, O.~E. Jensen, P.~K. Maini, and H.~M. Byrne.
\newblock A hybrid approach to multi-scale modelling of cancer.
\newblock {\em Philosophical Transactions of the Royal Society of London A:
  Mathematical, Physical and Engineering Sciences}, 368(1930):5013--5028, 2010.

\bibitem{Painter}
K.~J. Painter.
\newblock Continuous models for cell migration in tissues and applications to
  cell sorting via differential chemotaxis.
\newblock {\em Bulletin of Mathematical Biology}, 71:1117--1147, 2009.

\bibitem{Pao}
C.~V. Pao.
\newblock {\em Nonlinear parabolic and elliptic equations}.
\newblock Plenum Press, New York, 1992.

\bibitem{paulusal}
U.~Paulus, M.~Loeffler, J.~Zeidler, G.~Owen, and C.S. Potten.
\newblock The differentiation and lineage development of goblet cells in the
  murine small intestinal crypt: experimental and modelling studies.
\newblock {\em Journal of Cell Science}, 106:473--484, 1993.

\bibitem{Pretal}
S.~L. Preston, W.-M. Wong, A.~O.-O. Chan, R.~Poulsom, R.~Jeffery, R.A. Goodlad,
  N.~Mandir, G.~Elia, M.~Novelli, W.F. Bodmer, I.P. Tomlinson, and N.A. Wright.
\newblock Bottom-up histogenesis of colorectal adenomas: Origin in the
  monocryptal adenoma and initial expansion by crypt fission.
\newblock {\em Cancer Research}, 63:3819--3825, 2003.

\bibitem{RCS}
B.~Ribba, T.~Colin, and S.~Schnell.
\newblock A multiscale mathematical model of cancer, and its use in analyzing
  irradiation therapies.
\newblock {\em Theoretical Biology and Medical Modelling}, (3:7), 2006.

\bibitem{RMB}
L.~Roncucci, A.~Medline, and W.~R. Bruce.
\newblock Classification of aberrant crypt foci and microadenomas in human
  colon.
\newblock {\em Cancer Epidemiology, Biomarkers \& Prevention}, 1:57--60, 1991.

\bibitem{RCM2007}
T.~Roose, S.~J. Chapman, and P.~K. Maini.
\newblock Mathematical models of avascular tumor growth.
\newblock {\em SIAM Review}, 49(2):179--208, 2007.

\bibitem{Ross}
M.~H. Ross, G.~I. Kaye, and W.~Pawlina.
\newblock {\em Histology: A Text \& Atlas. Fourth Edition}.
\newblock Lippincott Williams \& Wilkins, 2003.

\bibitem{sancheznon}
E.~Sanchez-Palencia.
\newblock {\em Non-homogeneous Media and Vibration Theory}, volume 127.
\newblock Lecture Notes in Physics, Springer, Berlin, 1980.

\bibitem{SC}
J.~A. Sherratt and M.~A. Chaplain.
\newblock A new mathematical model for avascular tumour growth.
\newblock {\em Journal of Mathematical Biology}, 43:291--312, 2001.

\bibitem{Sal}
I-M. Shih et~al.
\newblock Top-down morphogenesis of colorectal tumors.
\newblock {\em Proceedings of the National Academy of Sciences of the United
  States}, 98(5):2640--2645, 2001.

\bibitem{SW}
B.~P. Stewart, B.~Wild, et~al.
\newblock {\em World Cancer Report}.
\newblock IARC Press, International Agency for Research on Cancer, 2014.

\bibitem{tartar2009general}
L.~Tartar.
\newblock The general theory of homogenization.
\newblock {\em Lecture Notes of the Unione Matematica Italiana}, 7, 2009.

\bibitem{tayloral}
R.~W. Taylor et~al.
\newblock Mitochondrial {DNA} mutations in human colonic crypt stem cells.
\newblock {\em The Journal of Clinical Investigation}, 112(9):1351--1360, 2003.

\bibitem{vLBJK}
I.~M.~M. van Leeuwen, H.~M. Byrne, O.~E. Jensen, and J.~R. King.
\newblock Crypt dynamics and colorectal cancer: advances in mathematical
  modelling.
\newblock {\em Cell Proliferation}, 39:157--181, 2006.

\bibitem{vLEIB}
I.~M.~M. van Leeuwen, C.~M. Edwards, M.~Ilyas, and H.~M. Byrne.
\newblock Towards a multiscale model of colorectal cancer.
\newblock {\em World Journal of Gastroenterology}, 13(9):1399--1407, 2007.

\bibitem{Van2009}
I.~M.~M. van Leeuwen et~al.
\newblock An integrative computational model for intestinal tissue renewal.
\newblock {\em Cell Proliferation}, 42(5):619--636, 2007.

\bibitem{walter}
A.~C. Walter.
\newblock {\em A Comparison of Continuum and Cell-based Models of Colorectal
  Cancer}.
\newblock PhD thesis, University of Nottingham, March 2009.

\bibitem{WK}
J.~P. Ward and J.~R. King.
\newblock Mathematical modelling of avascular-tumor growth.
\newblock {\em IMA Journal of Mathematics Applied in Medicine and Biology},
  14:39--69, 1997.

\end{thebibliography}

\end{document}